\documentclass[3p,prepritn]{elsarticle}

\usepackage[utf8]{inputenc}
\usepackage[T1]{fontenc} 
\usepackage{graphicx}%
\usepackage{multirow}%
\usepackage{amsmath,amssymb,amsfonts}%
\usepackage{amsthm}%
\numberwithin{equation}{section}
\numberwithin{figure}{section}
\usepackage[bb=dsserif,bbscaled=1.2]{mathalpha}
\usepackage{bm}
\usepackage{mathrsfs}%
\usepackage{mathtools}
\usepackage[title]{appendix}%
\usepackage{textcomp}%
\usepackage{manyfoot}%
\usepackage{tabularx}
\newcolumntype{C}{>{\centering\arraybackslash}X}
\newcolumntype{L}{>{\raggedright\arraybackslash}X}
\newcolumntype{R}{>{\raggedleft\arraybackslash}X}
\usepackage{booktabs}%
\usepackage{listings}%
\usepackage{todonotes}
\usepackage{pgfplots}
\usepackage{caption}
\usepackage{subcaption}
\usepackage{soul}
\usepackage{pgfgantt}
\usepackage{pifont}
\usepackage{csvsimple}
\usepackage{pgfplots}
\pgfplotsset{compat=1.18}
\usepackage[ruled,vlined,linesnumbered]{algorithm2e}
\usepackage{calc}
\usepackage{hyperref}
\usetikzlibrary{calc,positioning,patterns,shapes,backgrounds,fit,arrows.meta
}
\usepgfplotslibrary{fillbetween}
\usepgfplotslibrary{groupplots}
\usepgfplotslibrary{colorbrewer}
\pgfmathdeclarefunction{gauss}{2}{%
  \pgfmathparse{1/(#2*sqrt(2*pi))*exp(-((x-#1)^2)/(2*#2^2))}%
}
\newcommand{\pr}[2]{\textup{Pr}\left\{#1 \text{ in } \text{batch no. }#2\right\}}
\newcommand{\prp}[2]{\overline{\textup{Pr}}\left\{#1 \text{ in } \text{batch no. }#2\right\}}
\allowdisplaybreaks

\newcommand{\tikzxmark}{%
\tikz[scale=0.15] {
    \draw[line width=0.7,line cap=round] (0,0) to [bend left=6] (1,1);
    \draw[line width=0.7,line cap=round] (0.2,0.95) to [bend right=3] (0.8,0.05);
}}
\newcommand{\ppp}{\ensuremath{95^{\textrm{th}}}}

\pgfplotsset{colormap/greenyellow}

\begin{document}

\author[1]{Antonin Novak}
\ead{corresponding author: antonin.novak@cvut.cz}

\author[2]{Andrzej Gnatowski}
\ead{andrzej.gnatowski@pwr.edu.pl}

\author[1]{Premysl Sucha}
\ead{premysl.sucha@cvut.cz}

\affiliation[1]{Czech Institute of Informatics, Robotics and Cybernetics, Czech Technical University in Prague, CZ}

\affiliation[2]{Department of Control Systems and Mechatronics, Wroclaw University of Science and Technology, PL}

\title{Urgent Samples in Clinical Laboratories: Stochastic Batching to Minimize Patient Turnaround Time}

\begin{abstract}
    This paper addresses the problem of batching laboratory samples in hospital laboratories where samples of different priorities are received continuously with uncertain transportation times. The focus is on optimizing the control strategy for loading a centrifuge to minimize patient turnaround time~(TAT). While focusing on samples of patients in life-threatening situations (i.e., vital samples), we propose several online and offline methods, including a stochastic mixed-integer quadratic programming model integrated within a discrete-event system simulation. This paper aims to enhance patient care by providing timely laboratory results through improved batching strategies.
    The case study, which uses real data from a university hospital, demonstrates that incorporating distributional knowledge of transport times into our decision policy can reduce the median patient TAT of vital samples by 4.9 minutes and the 0.95 quantile by 9.7 minutes, but has no significant effect on low-priority samples.
    In addition, we show that this is essentially an optimal result by comparison with the upper bound obtained by a perfect-knowledge offline algorithm.
\end{abstract}

\maketitle

\section{Introduction}\label{sec:introduction}
Medical laboratories play an irreplaceable role in modern health care. 
Many diseases are diagnosed or confirmed by laboratory tests~\cite{Rohr2016-oo}. Moreover, early diagnosis is critical to prevent significant harm and even the death of a patient suffering from a disease such as an acute ischemic stroke or heart attack. Therefore, modern medical laboratories must be able to analyze samples quickly and provide results as soon as possible. 
An example of such a situation is when a patient with suspected myocardial infarction is brought into the hospital by an ambulance. The literature addressing such situations defines the door-to-balloon time \cite{Tsai2025} as the interval from the patient's arrival at the emergency department to the inflation of a balloon within the occluded coronary artery. This time should not exceed 90 minutes; otherwise, the risk of short-term mortality and major adverse cardiac events significantly increases. Myocardial infarction is diagnosed using the electrocardiogram. 
However, a troponin test is used to confirm or rule out myocardial damage if the electrocardiogram is inconclusive. In such a case, a blood sample is taken and immediately transported to the laboratory in vital mode, i.e., as a sample with the highest priority. In a situation where the team of specialists is waiting for confirmation or refutation of the diagnosis, every minute counts. For example, the authors in \cite{Scholz2018wx} reported that for patients with cardiogenic shock and out-of-hospital cardiac arrest, every 10-minute treatment delay resulted in 3.31 additional deaths among 100 percutaneous coronary intervention-treated patients.

The key performance indicator measuring how fast a laboratory can analyze a sample is called \textit{turn-around time} (TAT). 
Laboratories typically use \textit{laboratory TAT}, which is defined as the time between the moment when a sample is received by the laboratory and the time when the results are reported. Nevertheless, for patients or their physicians, this performance indicator is not relevant, as they are interested in draw-to-report TAT~\cite{Breil2011}, denoted in this paper as \textit{patient TAT}, which is defined as the time from collection of the sample from a patient to the reporting of the results.

The main reason why laboratories currently do not focus on patient TAT is that they cannot influence how samples are handled at individual wards and how quickly they are transported to the laboratory.
In other words, the collection and transport of samples to the laboratory are often not automated, while laboratories aim at total laboratory automation (TLA)~\cite{Genzen2018}, which minimizes human interaction with samples and thus indirectly makes the processing of samples much more predictable.
The laboratories are often equipped with fully automated systems, such as the Beckman Coulter DxA 5000, that can sort, centrifuge, and deliver samples to appropriate analyzers, perform prescribed tests, and report the results.
However, handling of samples is mostly not automated; even though information about sample collection is typically available in the laboratory information system, it is not yet exploited by laboratory automation systems to control patient TAT.

The research question we study in this paper is whether laboratory automation systems may benefit from information about samples that are on the way to the laboratory to minimize patient TAT. We focus on urgent samples that are usually associated with patients in life-threatening situations, e.g., acute myocardial infarction. 
Specifically, we study the processing of samples by centrifuge, which is carried out in batches, such as up to 56 samples in the case of the DxA 5000. 
If the centrifuge is started, it cannot be interrupted even if a priority sample has arrived or is about to arrive at the laboratory. A straightforward idea is to postpone the start of the centrifuge to process the incoming priority sample immediately. Nevertheless, the solution is not that easy in general. 
First, we cannot focus on a single sample, and more importantly, the transport time of samples to the laboratory is uncertain. Even if one knows that the sample was taken and is on the way, a decision requires consideration of many aspects, such as when the sample was taken, where it was taken, how it was transported (manually or by tube mail transport), and how likely different transport times are.

In this paper, we study an online sample batching problem, i.e., how to group clinical samples for centrifugation and when the centrifuge should be started. 
In this research, we assume that there is a single centrifuge. 
We are inspired by the laboratory automation system DxA~5000, which typically features a single centrifuge, ideal for small and medium-sized hospitals. 
Another reason is the space limitations in laboratories, which typically do not allow the installation of multiple configuration modules, which would result in a higher system cost.
To minimize patient TAT, we exploit online information given when a sample is collected from a patient. The transport time of a sample is modeled as a random variable characterized by a probability density function. The problem is formulated as a stochastic optimization problem minimizing patient TAT for the samples with the highest priority, denoted as \textit{vital}; this category is reserved for patients in immediate danger of death. Lower priorities are \textit{statim}, which is used for less urgent situations than vital situations, e.g., suspicion of a disease, and \textit{routine}, for ordinary samples, where the result can typically be delivered the next day. The benefit of our approach is demonstrated using real data from University Hospital Královské Vinohrady, Czech Republic.
We show that a simple rule-based solution deteriorates patient TAT for high-priority samples, whereas our solution provides significantly better solutions and delivers more predictable (consistent) patient TAT.
\subsection{Contributions and outline}
In this paper, we aim to increase the efficiency and reliability of TLA systems with mathematical optimization tools to improve the dispatching of the centrifuge.
Specifically, the main contributions are as follows:
\begin{enumerate}
    \item We propose a new way vital samples can be processed in the classic laboratory automation systems supporting only statim and routine sample priorities;
    \item We develop optimization techniques that focus on the patient TAT rather than laboratory TAT, increasing the impact of the decision directly for the patient;
    \item We propose an online stochastic mathematical model that uses information about transiting samples to optimally control batching and the start times of the centrifuge;
    \item We validate our approaches using discrete event simulation on a case study that uses real data obtained from a university hospital;
    \item The experiments show that our online dispatching policy improves the 0.95 quantile of patient TAT compared to the existing solution by 9.7 minutes by using distributional knowledge of the uncertain transport times combined with the stochastic mathematical optimization model;
    \item The proposed perfect-knowledge offline algorithm demonstrates that the performance of the online dispatching policy is essentially optimal;
    \item We provide a complexity characterization of the problem in the offline setting.
\end{enumerate}
The paper is structured as follows. 
In Section~\ref{sec:relwork}, we survey related works concerning the problem considered and the optimization techniques used.
The formal problem statement is given in Section~\ref{sec:problem_statement}.
Our dispatching policies are proposed in Section~\ref{sec:methods} and are benchmarked using our discrete event simulation system (Section~\ref{sec:des}) in Section~\ref{sec:experiments}.

In the text below, we use the following notation.
The uncertain parameters that are treated as random variables are denoted with the tilde symbol $\tilde{x}$, where their realizations are denoted with the hat symbol $\hat{x}$.
The other parameters that are assumed to be deterministic are denoted in plain $x$. 
See Table~\ref{tab:notation} for a detailed description of the parameters used.
\section{Related work}\label{sec:relwork}
The review of the related literature is split into two sections to reflect two perspectives: (i)~the application point of view, reflecting the work on the modeling and optimization of laboratory automation, including its impact on the TAT, and (ii)~related combinatorial problems from the area of batching, simulation and optimization paradigms, scheduling with uncertain release and transport times and online (dynamic) scheduling problems.

\subsection{Optimization of Laboratories}
In the area of the optimization of automated laboratory systems, it is essential to define key performance indicators, which may differ substantially among the various stakeholders.
For example, commercial laboratories that focus on routine tests primarily minimize operating costs (e.g., material consumption, method duplication), whereas a hospital laboratory typically minimizes the TAT of priority (vital and statim) samples because the lives of the patients depend on these.
However, as mentioned in Section~\ref{sec:introduction}, laboratories currently keep track only of the laboratory TAT, whereas it was shown that the preanalytical phase is the most error-prone~\cite{tsai2019critical}; thus, it should not be excluded from the TAT calculation.
In our work, we follow this line of reasoning; thus, we focus on minimizing the patient TAT instead of the laboratory TAT. 
A crucial consideration that is equally important is which function of patient TAT to minimize.
For a set of processed samples, we observe an empirical distribution of patient TATs---the time to complete the results is not deterministic, as it depends on the transportation time, waiting time in the buffer, congestion at analyzers, etc.
Again, different parties may opt to optimize the expected patient TAT, the worst-case TAT, or generally an arbitrary quantile of the TAT distribution. 
In \cite{holland2005reducing}, the impact of the TAT on patient length of stay was analyzed.
It was shown that the mean TAT does not well correspond to the length of stay but rather is affected by the so-called TAT outlier percentage, which reflects the number of samples that were completed well beyond the nominal level.
This finding suggests that reducing patient TAT at high quantile levels contributes to shorter patient length of stay.

A problem similar to our setting was addressed by
Yang et al.~\cite{yang2015optimization}, where the authors identified limitations and bottlenecks of laboratories in a medical center in Taiwan. 
The analysis revealed high waiting times for the three identical DxC analyzers, which were caused by an uneven distribution of samples among the analyzers. 
They optimized the controls for the batch waiting time of the centrifuge (i.e., how long the centrifuge waits for a whole batch) and the distribution of samples among the DxC analyzers. 
The typical classes of samples are distributed among DxC analyzers on the basis of their average processing time with respect to two defined value thresholds. 
Mixed-integer programming (MIP) is utilized to find optimal values for thresholds and the waiting time of the centrifuge to maximize the number of samples with TATs less than 60 minutes, yielding a 54\% improvement in the average sample TAT. 
 This work disregards the priorities of samples and transportation times; nevertheless, the main distinction between their work and ours is that they do not consider the centrifuge dispatching algorithm in an online setting. Instead, they derive optimal timeout thresholds from historical data and then set them statically.
On the other hand, we treat the centrifuge control in an online setting and adjust the decision optimally with respect to the available information about transiting samples.
\begin{figure}
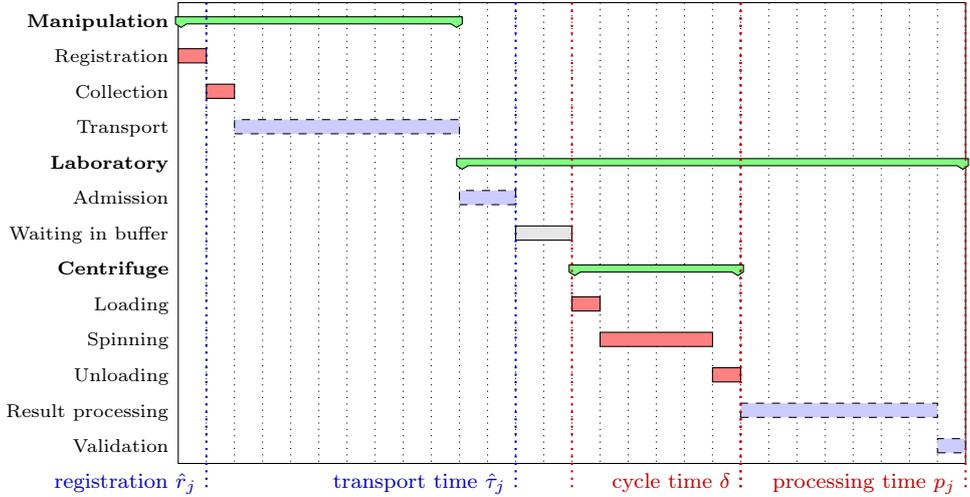

    \centering
    \begin{ganttchart}[
       vgrid,
    bar/.append style={fill=red!50},
    group/.append style={draw=black, fill=green!50},
    bar label anchor/.append style={align=left, text width=30em},
    group label anchor/.append style={align=left, text width=30em},
    y unit chart=0.47cm,
    x unit=0.37cm,
     transportlabel/.style={
        vrule label node/.append style={anchor=north east},
        vrule label font={\footnotesize},
        vrule/.style={dotted, thick, blue!75!black, shorten >=-1em},
    },
     spinninglabel/.style={
        vrule label node/.append style={anchor=north east},
        vrule label font={\footnotesize},
        vrule/.style={dotted, thick, red!75!black, shorten >=-1em},
    },
     processinglabel/.style={
        vrule label node/.append style={anchor=north east},
        vrule label font={\footnotesize},
        vrule/.style={dotted, thick, cyan!75!black, shorten >=-1em},
    }
    ]{1}{28}
        \ganttgroup{\scriptsize Manipulation}{1}{10} \\
        \ganttbar{\scriptsize Registration}{1}{1} \\
        \ganttbar{\scriptsize Collection}{2}{2} \\
        \ganttbar[bar/.append style={fill=blue!20, dashed}]{\scriptsize Transport}{3}{10} \\
        \ganttgroup{\scriptsize Laboratory}{11}{28} \\
        \ganttbar[bar/.append style={fill=blue!20, dashed}]{\scriptsize Admission}{11}{12} \\
        \ganttbar[bar/.append style={fill=gray!20}]{\scriptsize Waiting in buffer}{13}{14} \\
        \ganttgroup{\scriptsize Centrifuge}{15}{20} \\
         \ganttbar{\scriptsize Loading}{15}{15}\\
         \ganttbar{\scriptsize Spinning}{16}{19}\\
         \ganttbar{\scriptsize Unloading}{20}{20}\\
        \ganttbar[bar/.append style={fill=blue!20, dashed}]{\scriptsize Result processing}{21}{27}\\
        \ganttbar[bar/.append style={fill=blue!20, dashed}]{\scriptsize Validation}{28}{28}
        \ganttvrule[transportlabel]{registration $\hat{r}_j$}{1}
        \ganttvrule[transportlabel]{\footnotesize transport time $\hat{\tau}_j$}{12}
        \ganttvrule[spinninglabel]{}{14}
        \ganttvrule[spinninglabel]{\footnotesize cycle time $\delta$}{20}
        \ganttvrule[spinninglabel]{}{20}
        \ganttvrule[spinninglabel]{\footnotesize processing time $p_j$}{28}
    \end{ganttchart}
    \caption{Timing diagram of the laboratory process for a patient sample. The red rectangles correspond to activities whose duration is essentially deterministic, whereas the blue rectangles correspond to activities whose duration is uncertain. The light gray depicts the idle time of a sample.}
    \label{fig:process-timing}
\end{figure}

Similarly, Sojma~\cite{Sojma_2022} utilizes a discrete event system (DES) to describe the workflow of the DxA 5000 system, encompassing the complete sample life cycle from input to centrifuging, sample transport on the conveyor belt, and method testing at various analyzers.
The optimization aimed to redistribute the allocation of methods among the laboratory analyzers, balancing their utilization and increasing the system throughput.
For this purpose, a simplified surrogate MIP formulation of the system was proposed, with the DES being used to validate the newly proposed method assignment. 
The simulation confirmed that the new assignment improved the throughput of biochemical analyzers by 11.7\%.

Laboratory processes share certain characteristics with the production domain. Thus, similar optimization approaches, objective functions, and modeling paradigms can be used.
For example, \cite{li2019medical} represents a laboratory workflow as a production process with the optimization of so-called weighted flow time minimization.
They consider each analyzer as a machine and the patient's sample as a job to be processed on a machine.
They designed an online policy that batches the samples of different priorities to be processed by a set of laboratory analyzers to minimize the maximum weighted flow time. 
The online policy essentially achieves optimal performance among all online policies for the case in which the batch has infinite capacity.
Under the bounded capacity model, their policy has optimal performance for weights $\Omega_j\in \left[1,2\right]$.
Work~\cite{li2019medical} is particularly relevant to our problem. However, the main difference is that in our work, we focus on the preanalytic stage of the process (i.e., centrifuge rather than biochemical analyzers). In addition, the model presented in their work disregards the transport times of samples, and the processing time of samples is always one.

With respect to algorithmic tools, MIP is a successful technique for optimizing health care processes. 
For example, Farhadi et al. \cite{farhadi2023optimization} used MIP for scheduling patient appointments with technicians in hemodialysis centers.
In cases where some parameters of the problem are uncertain, theoretical frameworks, such as stochastic programming or distributional robust optimization, must be deployed for uncertainty mitigation ~\cite{NOVAK2022438}.
For example, to address demand uncertainty, \cite{lalmazloumian2023two} employed a two-stage stochastic programming framework to reformulate a stochastic formulation of operating room allocation.
Alternatively, predictive models~\cite{schafer2023combining,shi2023surgical} can be used to estimate the parameter distribution and employ those in a data-driven mathematical programming model.
However, some health care management processes are too complex to be modeled by an MIP model. 
In these cases, the processes can often be modeled by a DES~\cite{osorio2017simulation} to provide high-fidelity simulation.
However, the direct optimization of the parameters within such a system often becomes intractable; thus, the approaches often resort to using DES simulation as a fitness function oracle in a black-box optimization or using DES as a tool to verify the optimization results from a surrogate optimization problem that can be solved.
For example, Osorio et al.~\cite{osorio2017simulation} employed the iterative optimization-based simulation scheme to solve the production planning problem in the blood supply chain, where they iteratively used DES to compute values that act as the inputs to the MIP for the following iterations.
\begin{table}
    \centering
    \caption{Overview of the used parameters and symbols.}
    \label{tab:notation}
    \begin{tabular}{c|c|c}
    \toprule
         parameter & description & value \\
         \midrule
          $\tilde{r}_j$ &  sample registration time & uncertain \\
          $\hat{r}_j$ & realization of registration time & deterministic\\
            $\tilde{\tau}_j$ &  sample transportation time & uncertain\\
          $\hat{\tau}_j$ & realization of transportation time & deterministic\\
          $\tilde{a}_j$ & sample arrival time & uncertain \\
          $\hat{a}_j$ &realization of arrival time & deterministic \\
          $\lambda_j$ & sample priority &deterministic\\
          $w_j$ & hospital ward of sample registration & deterministic\\
          $p_j$ & sample processing time & deterministic \\
          $C_j$ & sample completion time & deterministic \\
          $S^B$ & batch start time & deterministic \\
          $\delta$ & centrifuge cycle time & deterministic \\
          $V_t$ & set of vitals available at time $t$ & deterministic\\
          $\overline{V}_t$ & set of vitals in transit at time $t$ & deterministic\\
         \bottomrule
    \end{tabular}
    
\end{table}
\subsection{Batching and allocation problems in stochastic and online environments}
Timing the centrifuge runs resembles a parallel batching (p-batching) problem~\cite{fowler2022survey}.
In a typical setting, the jobs have processing times, release times, and due dates, and the goal is to group available jobs into a batch with a limited capacity, whose processing time is given as the maximum of the processing times of jobs within the batch.
The goal is often to minimize a function of the completion times of the jobs.
In that respect, perhaps the most similar problem to ours is the setting where we consider a single machine (centrifuge), jobs with constant processing times (constant spinning duration of the centrifuge), release times, and due dates.
This problem under an offline, deterministic setting was studied by Baptiste~\cite{Baptiste2000}, who reported that the problem is solvable in polynomial time for various natural objective functions.
Later, his result was improved by~\cite{Condotta2010}, who proposed a much faster and thus more practical algorithm for the feasibility variant of the problem.
To support comparison with the perfect knowledge offline algorithm, we build on their work.
However, the main differences in optimizing the laboratory centrifuge from the problems described above are the presence of uncertain parameters, such as release and transport times.

A recent survey on offline p-batch problems was conducted by Fowler and M\"{o}nch~\cite{fowler2022survey}.
As one of their conclusions, they argue for the increased need to study batching problems in a stochastic setting.
Indeed, as noted by~\cite{rocholl2022scheduling} \emph{``there are only a very few papers that deal with the p-batch scheduling problem involving uncertain data''}.
Therefore, \cite{rocholl2022scheduling} consider a p-batching problem with incompatible jobs (i.e., some families of jobs cannot be processed in the same batch) and uncertain ready times with the total weighted tardiness criterion.
They introduce sampling-based heuristics that assume a distribution of delays in the release times to create a robust baseline schedule.
However, owing to the unforeseen realizations of the release times, the schedule must be repaired through an online policy as the realizations unfold.

In fact, our problem also resembles a variant of a stochastic p-batching problem with uncertain release times.
Nevertheless, the difference is that the decision-maker learns about sample $j$ at the time it is registered in the system $\hat{r}_j$, but it becomes available for centrifuging when it arrives at the laboratory after an uncertain duration $\tilde{\tau}_j$; thus, the release time is $\tilde{r}_j = \hat{r}_j+\tilde{\tau}_j$.
Moreover, the registration times $\hat{r}_j$ are revealed gradually, without a particular assumption regarding their distribution.
Therefore, our work is also related to the stream of works considering the online setting of the batching problem~\cite{CALS2021107221}, where the set of jobs is not known in advance but is released over time~\cite{XIONG2024122752, li2019medical}.
The goal is to design a policy that makes online (here-and-now) decisions about which jobs should be batched together and when to start the batch.
This problem frequently arises in semiconductor manufacturing, wafer fabrication~\cite{KOO20131690}, and e-commerce order fulfillment~\cite{pardo2024order}.
Koo and Moon~\cite{KOO20131690} classify the dispatching policies into two main categories: (i)~threshold-based policies that are applied in situations where no information except the observable state of the system is available and (ii)~look-ahead policies that use information on the near-future system state.
In this paper, we describe various centrifuge policies with increasing amounts of utilized information, exactly matching the above classification of the dispatching rules.
Additionally, the problem of dispatching the centrifuge can be considered through the lens of the control theory of discrete systems, which has been recently addressed by the reinforcement learning paradigm~\cite{CALS2021107221}.
In this context, the task is to derive a controller represented by a neural network model through extensive simulation of the system and the environment.
The model observes samples from the observation space (e.g., sensor values) and performs certain actions with its actuators (e.g., required motor torque). 
However, in our setting, the observation space is not fixed in size, as there might be an arbitrary number of transiting samples at a given time.
Reinforcement learning models and algorithms typically disregard this aspect.
Furthermore, even though it was empirically observed that such reinforcement learning controllers might achieve good performance in the expected sense, how to verify the trained policies they represent~\cite{li2023sok} to obtain guarantees for the worst-case scenario remains an open problem.

Some aspects of the problem, such as uncertain arrival time, also appear in queuing theory~\cite{harchol2021open,sulz2024potential}, where one of the targets is to directly minimize response tail time or the $q$-th quantile of response time, i.e., the time between the completion and the release of the job.
The problem of sharing system resources (i.e., centrifuges) among activities with different criticalities (i.e., samples with different priorities) has been studied in mixed-criticality systems~\cite{NOVAK2019687}.
 
\section{Problem statement} \label{sec:problem_statement}
We assume a set of $n$ samples $\mathcal{S}=\{1, 2, \ldots, j, \ldots, n\}$ to be released over the horizon of 24 hours.
Each sample $j\in\mathcal{S}$ is characterized by its registration time $\tilde{r}_j$, transportation time $\tilde{\tau}_j$, and sample priority $\lambda_j \in \{\textsc{routine},\ \textsc{statim},\ \textsc{vital}\}$ with \textsc{routine} being the least and \textsc{vital} being the most important priority.
Furthermore, the hospital ward where the sample is taken is denoted as $w_j\in\{\textsc{ward}_1, \textsc{ward}_2, \ldots, \textsc{ward}_W\}$, where $W$ represents the total number of wards in a hospital. 
Finally, sample $j$ is associated with processing time $p_j$, which represents the time the laboratory system needs to provide the result from the centrifuged sample.
The number of $n$ samples, their registration times $\tilde{r}_j$, and their transport times $\tilde{\tau}_j$ are random variables; thus, the realizations are not known beforehand. 
When the sample is registered in the system, we observe the realization of its registration (release) time $\hat{r}_j$. 
Furthermore, we assume that by registering sample $j$ into the system, we also learn the realization of sample priority $\lambda_j$, the hospital ward $w_j$, and its processing time $p_j$.
The realization of transportation time $\hat{\tau}_j$ is not revealed until sample~$j$ arrives at the laboratory.
To simplify the notation, we also introduce the sample arrival time $\hat{a}_j = \hat{r}_j + \hat{\tau}_j$.

Every sample $j\in \mathcal{S}$ needs to be centrifuged before it is processed by laboratory analyzers.
The processing time of a sample depends on the methods that should be performed on the sample.
We assume that when the sample is registered, we learn the list of methods that can be used to estimate its processing time $p_j$ in the analytical part of the system.
The capacity of the centrifuge is $N=56$, and its cycle time, including spinning time, sample loading, and sample unloading, is $\delta=15\cdot 60$ seconds.
When the centrifuge is started, its run cannot be interrupted, and it must complete its cycle. 
Therefore, when sample $j$ is included in the current batch and the centrifuge starts its cycle at time $t$, the sample is processed at time $C_j=t + \delta + p_j$.
We call patient turn-around time (patient TAT) the time between sample completion and registration (sample collection); i.e., the turn-around time of sample $j$ is $\text{TAT}_j = C_j - \hat{r}_j$.
Thus, the patient TAT is the sum of the manipulation time plus the laboratory TAT.
Note also that the notion of patient TAT directly corresponds to the term flow time $F_j$ typically used in production scheduling.
The main symbols and notations used are summarized in Table~\ref{tab:notation}.

\begin{figure}
    \centering
    \includegraphics[width=0.57\textwidth]{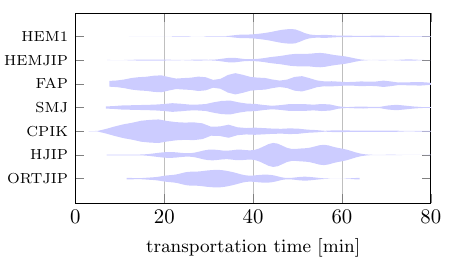}
    \caption{Transport time probability densities estimated by kernel density estimation of statim samples for different hospital wards.}
    \label{fig:transport-distribution}
\end{figure}

An example of a sample life cycle is illustrated in Figure~\ref{fig:process-timing}.
The red rectangles correspond to activities that have essentially deterministic durations.
The blue rectangles represent activities with uncertain duration, i.e., transportation and processing time.
Here, we model the sample transportation time $\tilde{\tau}_j$ using a probability distribution, which is assumed to be a function of sample priority $\lambda_j$ and its hospital ward $w_j$.
For example, see the transport time distribution for $\lambda_j=\textsc{statim}$ for seven different wards at University Hospital Královské Vinohrady, depicted in Figure~\ref{fig:transport-distribution}.
For simplicity, we further assume that the duration of registration and sample collection is included within $\tilde{\tau}_j$.
The processing time $p_j$ of the sample is influenced primarily by the list of methods to be performed and the utilization of the analyzers, as well as their potential outages.
In this paper, we do not focus on modeling the processing component of the system.
Instead, we use a coarse approximation that estimates $p_j$ as the processing time of the longest method to be performed on sample~$j$ (i.e., the requested methods of the sample are evaluated in parallel).
The light gray rectangle in Figure~\ref{fig:process-timing} depicts the idle time of the sample and is subject to optimization.
In this paper, we aim to reduce the sample idle time by suitable batching decisions and centrifuge timing so that it decreases the patient TAT.

The main key performance indicator of the laboratory is the expected TAT of vital samples to ensure timely delivery of the results.
However, the laboratories also keep track of the 0.95 TAT quantile of vital samples to prevent too many samples from being delayed in processing.
The value 0.95 is typically chosen to discard outliers from the statistical indication that occasionally occur in the analytical part of the process (e.g., reruns of inconclusive test results).
Next, in order of importance, the expected and 0.95 quantile of the statim sample is considered.
To model the problem, we use the following assumptions.
We assume that the sample transportation time is a function of the priority and the hospital ward, 
i.e., $\tilde{\tau}_j \sim f(\lambda_j, w_j)$.

However, in principle, our approach works with any model that estimates the probability density function of $\tilde{\tau}_j$~\cite{shi2023surgical} or provides a single-point estimate with confidence intervals.
Without loss of generality, we assume that the samples $\mathcal{S}$ are about to be released (registered) over the 24-hour time horizon.
Finally, we assume that when the centrifuge starts loading, the set of samples to be loaded is fixed and cannot change.
Therefore, it can be seen as a process where all samples in the batch are loaded/unloaded at the same time, taking a constant time.

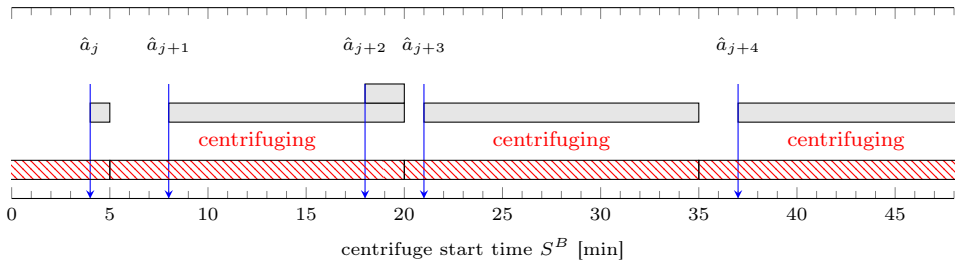
\begin{figure}
    \centering
    \begin{tikzpicture}
        \begin{axis}[
            height=0.25\textwidth,
            width=0.85\textwidth,
            ymin=0,
            ymax=0.5,
            xmin=0,
            xmax=48,
            xlabel={centrifuge start time $S^B$ [min]},
            label style={font=\scriptsize},
            tick label style={font=\scriptsize},
            minor x tick num=4,
            hide y axis
            ]
            \draw[fill=gray!20] (axis cs:37,0.2) rectangle (axis cs:50,0.25);
            \draw[fill=gray!20] (axis cs:4,0.2) rectangle (axis cs:5,0.25);
            \draw[fill=gray!20] (axis cs:8,0.2) rectangle (axis cs:20,0.25);
            \draw[fill=gray!20] (axis cs:18,0.25) rectangle (axis cs:20,0.30);
            \draw[fill=gray!20] (axis cs:21,0.2) rectangle (axis cs:35,0.25);
            \draw [->,>=stealth,blue] (axis cs:4,0.3) -- (axis cs:4,0.0);
            \node[rotate=0] at (axis cs:4,0.4){\scriptsize{$\hat{a}_j$}};
            
            \draw [->,>=stealth,blue] (axis cs:8,0.3) -- (axis cs:8,0.0);
            \node[rotate=0] at (axis cs:8,0.4){\scriptsize{$\hat{a}_{j+1}$}};
            \draw [->,>=stealth,blue] (axis cs:18,0.3) -- (axis cs:18,0.0);
            \node[rotate=0] at (axis cs:18,0.4){\scriptsize{$\hat{a}_{j+2}$}};
            \draw [->,>=stealth,blue] (axis cs:21,0.3) -- (axis cs:21,0.0);
            \node[rotate=0] at (axis cs:21,0.4){\scriptsize{$\hat{a}_{j+3}$}};
            \draw [->,>=stealth,blue] (axis cs:37,0.3) -- (axis cs:37,0.0);
            \node[rotate=0] at (axis cs:37,0.4){\scriptsize{$\hat{a}_{j+4}$}};
            \draw[pattern=north west lines, pattern color=red] (axis cs:-10,0.05) rectangle (axis cs:5,0.1);
            \draw[pattern=north west lines, pattern color=red] (axis cs:5,0.05) rectangle (axis cs:20,0.1);
            \node[red] at (axis cs:12.5,0.15){\footnotesize{centrifuging}};
            \draw[pattern=north west lines, pattern color=red] (axis cs:20,0.05) rectangle (axis cs:35,0.1);
            \node[red] at (axis cs:27.5,0.15){\footnotesize{centrifuging}};
            \draw[pattern=north west lines, pattern color=red] (axis cs:35,0.05) rectangle (axis cs:50,0.1);
            \node[red] at (axis cs:42.5,0.15){\footnotesize{centrifuging}};
        \end{axis}
    \end{tikzpicture}
    \caption{A fixed-schedule centrifuge policy that is not triggered by sample arrival may miss an arriving sample. The idle times of the samples are depicted in gray.}
    \label{fig:ex_fixsched}
\end{figure}

One of the simplest solutions to the above problem would be to start the centrifuge every $\delta$ seconds, thus achieving its highest utilization. 
However, such a strategy, which we refer to as the \textit{fixed-schedule policy}, can be inefficient from the perspective of the patient TAT.
Indeed, consider the example in Figure~\ref{fig:ex_fixsched}, where blue arrows correspond to time instances where a sample enters the laboratory, and the hatched red parts correspond to the intervals where the centrifuge is occupied.
If the timing of the centrifuge is not optimized, the sample may need to wait in the buffer to be processed in the next centrifuge run, as shown by the gray rectangles.
In this paper, we aim to design control strategies that mitigate the long waiting times of the samples with the goal of minimizing the patient TAT of samples with respect to the key performance indicators mentioned above.
\begin{table}
    \centering
    \caption{Overview of the analyzed and proposed methods.}
    \label{tab:alg_list}
    \begin{tabular}{c|ccc}
        \toprule
         & \multicolumn{3}{c}{environment information}\\
         & available & registered & transport \\
        method & samples & samples & distribution \\
        \midrule
          fixed-schedule policy, Section~\ref{sec:problem_statement} & \tikzxmark &  \tikzxmark & \tikzxmark \\
         threshold-based policy, Section~\ref{sec:sota_alg}  & \checkmark & \tikzxmark &  \tikzxmark\\
         look-ahead policy, Section~\ref{sec:fifo_alg} & \checkmark & \checkmark & \tikzxmark \\
         stochastic MIQP policy, Section~\ref{sec:mip_alg} & \checkmark  & \checkmark & \checkmark \\
         offline perfect-knowledge, Section~\ref{sec:perfect_knowledge} & N/A  & N/A & N/A \\
         \bottomrule
    \end{tabular}
\end{table}

\section{Batching algorithms}\label{sec:methods}
In this section, we present three batching methods that address the problem of optimal centrifuge timing and sample batching.
Their high-level strategies and capabilities differ mainly in the amount of information utilized by the algorithm, ranging from observing only the samples waiting inside the laboratory to assuming the distribution of the transport times of the samples.
Table~\ref{tab:alg_list} summarizes the algorithms proposed in this work.
Additionally, the distinction between the methods that are agnostic to transporting samples and those that utilize transport time distribution can be seen as the difference between optimizing laboratory TAT versus patient TAT.

Section~\ref{sec:des} describes the \textit{discrete event system} simulation, which represents a unified interaction of the batching policies with the process outlined in Figure~\ref{fig:process-timing}. The following subsections address the batching methods.
Arguably, the simplest imaginable batching method is the fixed-schedule policy demonstrated in Figure~\ref{fig:ex_fixsched}, which runs the centrifuge every $\delta$ seconds.
An improved solution is presented in Section~\ref{sec:sota_alg}, where we describe an algorithm that is often used in laboratories and is considered the baseline in this paper.
According to the batching policy classification, e.g., by \cite{KOO20131690, yang2015optimization}, when no information except the current system state is available, the policy is referred to as a \textit{threshold-based policy}.
In Section~\ref{sec:fifo_alg}, we introduce our algorithm, which addresses some of the disadvantages of the baseline solution.
The algorithms that make use of the near-future state of the system are referred to as \textit{look-ahead policies}.
Finally, in Section~\ref{sec:mip_alg}, we propose a stochastic mathematical optimization model that uses the transport time distribution to minimize the patient TAT of vital samples.
To derive the best theoretical performance achievable, in Section~\ref{sec:perfect_knowledge}, we develop an offline perfect-knowledge algorithm that has access to all future realizations of registration and transport times.
\subsection{Discrete event system simulation}\label{sec:des}
To verify the functionality and evaluate the efficiency of the proposed methods, we designed a discrete event simulation~\cite{osorio2017simulation}. 
It provides an environment resembling the process outlined in Figure~\ref{fig:process-timing}, i.e., the registration of new samples into the system, their transport, the behavior of the centrifuge with its control algorithm, and the calculation of the patient TAT. 
Therefore, all the proposed and investigated control algorithms for the centrifuge follow an identical protocol for interaction with the environment and for evaluation.

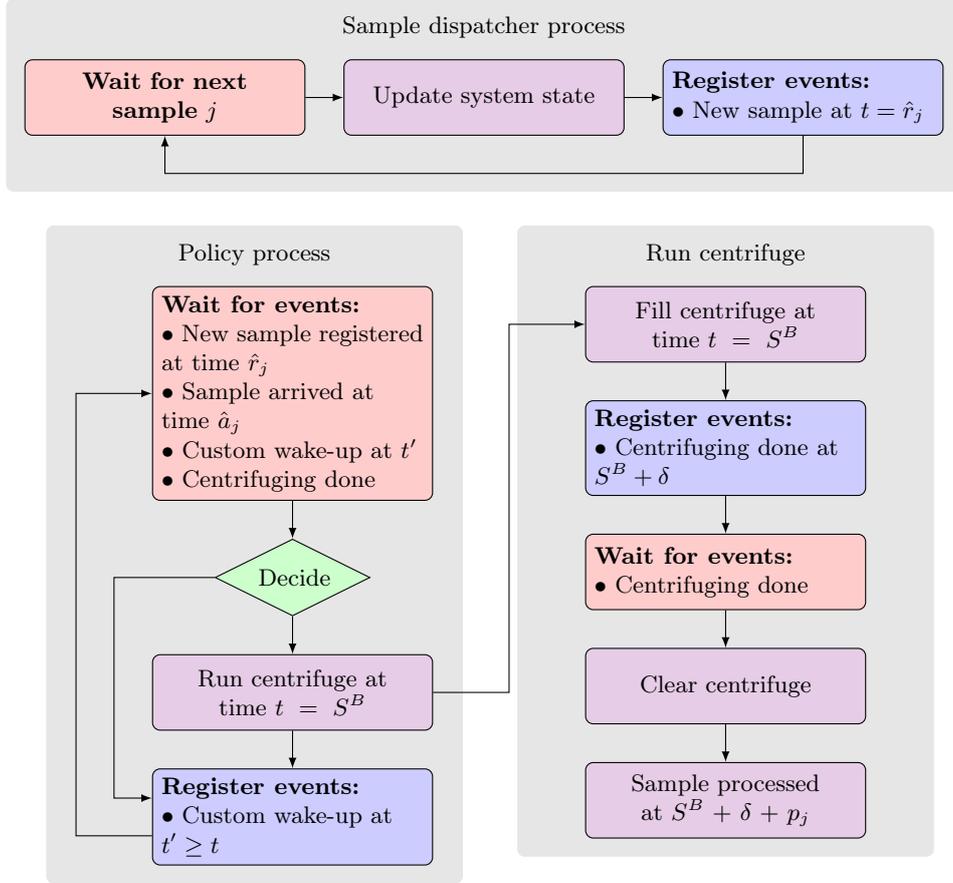
\begin{figure}[t]
    \centering
    \newcommand{\around}[2]{%
    \node[fit=#1] (b1) {};
    \node[above=1pt of b1, align=center, minimum height=0cm, text width=5cm] (bl1) {#2};
    \node[fit=(b1)(bl1), bkg] {};
    \node[above=1pt of b1, align=center, minimum height=0cm, text width=5cm] (bl1) {#2};
}
\newcommand{\nmh}{1cm}
\begin{tikzpicture}[
        node distance=.5cm and .5cm,
        font=\small,
        every node/.style={%
            align=flush left, 
            anchor=north,
            text width=3.45cm, 
            minimum height=\nmh, 
            rounded corners=0.1cm
        },
        wait/.style={draw,fill=red!20},
        stateupdate/.style={draw,fill=violet!20,align=center},
        registerevent/.style={draw,fill=blue!20},
        other/.style={draw,align=center, fill=white},
        startstop/.style={rectangle,draw,align=center,fill=green!20},
        decide/.style={
            diamond, 
            draw,
            fill=green!20,
            text width=0.9cm,
            align=center, 
            minimum height=\nmh, 
            aspect=2,
            rounded corners=0cm
        },
        bkg/.style={rectangle, fill=black!10},
        normalarc/.style={-latex},
    ]
        
        \node[wait,align=center] (l21) at (0,0) {\textbf{Wait for next\\ sample $j$}};
        \node[stateupdate, right= of l21] (l22) {Update system state};
        \node[registerevent, right= of l22] (l23) {%
            \textbf{Register events:}%
            \\$\bullet$ New sample at~$t=\hat{r}_j$};
        \draw[normalarc] (l21) -- (l22);
        \draw[normalarc] (l22) -- (l23);
        \coordinate (a1) at ($(l23.south)!0.5!(l21.south)+(0,-.5cm)$);
        \draw[normalarc] (l23) |- (a1) -| (l21);
        
        \node[wait, below=2.5 cm of $(l21)!.4!(l22)$] (l31) {%
            \textbf{Wait for events:}%
            \\$\bullet$ New sample registered at time $\hat{r}_j$
            \\$\bullet$ Sample arrived at time $\hat{a}_j$
            \\$\bullet$ Custom wake-up at $t^\prime$
            \\$\bullet$ Centrifuging done};
        \node[decide, below= of l31] (l32) {Decide};
        \node[stateupdate, below= of l32] (l33) {Run centrifuge at time $t=S^B$};
        \node[registerevent, below= of l33] (l34) {%
            \textbf{Register events:}%
            \\$\bullet$ Custom wake-up at $t^\prime\geq t$};
        \draw[normalarc] (l31) -- (l32);
        \draw[normalarc] (l32) -- (l33);
        \draw[normalarc] (l33) -- (l34);
        \coordinate (a2) at ($(l31.west)!0.5!(l34.west)+(-1cm,0)$);
        \coordinate (a3) at ($(l31.west)!0.5!(l34.west)+(-.5cm,0)$);
        \draw[normalarc] ($(l34.west)+(0,-0.25cm)$) -| (a2) |- (l31);
        \draw[normalarc] (l32.west) -| (a3) |- ($(l34.west)+(0,0.25cm)$);
        \node[stateupdate, right= 2cm of l31.north east, anchor=north west] (l41) {Fill centrifuge at time $t=S^B$};
        \node[registerevent, below= of l41] (l42) {%
            \textbf{Register events:}%
            \\$\bullet$ Centrifuging done at $S^B+\delta$};
        \node[wait, below= of l42] (l43) {%
            \textbf{Wait for events:}%
            \\$\bullet$ Centrifuging done};
        \node[stateupdate, below= of l43] (l44) {Clear centrifuge};
        \node[stateupdate, below= of l44] (l45) {Sample processed at $S^B+\delta+p_j$};
        \draw[normalarc] (l41) -- (l42);
        \draw[normalarc] (l42) -- (l43);
        \draw[normalarc] (l43) -- (l44);
        \draw[normalarc] (l44) -- (l45);
        \draw[normalarc] (l33) -| ($(l33)!0.5!(l41)$) |- (l41);
        \begin{pgfonlayer}{background}
            \around{(l21)(l22)(l23)(a1)}{Sample dispatcher process}
            \around{(l31)(l32)(l33)(l34)(a2)(a3)}{Policy process}
            \around{(l41)(l42)(l43)(l44)(l45)}{Run centrifuge}
        \end{pgfonlayer}
    \end{tikzpicture}
    \caption{Discrete event simulation modeling the sample transit and centrifuging in a clinical laboratory.}
    \label{fig:DES}
\end{figure}

A discrete event system (DES) is a model in which the state evolution is determined by the occurrence of a finite number of events over time. 
The time is discrete, and between events, the system state cannot change.
To simulate processing samples in a clinical laboratory, we designed a DES, as summarized in the diagram in Figure~\ref{fig:DES}. 
It consists of two processes--the sample dispatcher process and the policy process. 
The policy process controls the centrifuge and is invoked either when (i)~a new sample is registered, (ii)~a sample arrives, (iii)~a centrifuging cycle is completed, or (iv)~at a custom time $t^\prime$, if requested by the policy process.
\subsection{Threshold-based policy} \label{sec:sota_alg}
The baseline algorithm that is often deployed in laboratories~\cite{yang2015optimization} can be classified as a priority threshold-based batching algorithm.
The threshold-based policy consists of a set of rules that decide whether to postpone or dispatch the centrifuge on the basis of the information about the samples that are already available in the laboratory.
The advantage of rule-based methods is their explainability and predictable behavior, but as will be demonstrated later, they may lead to suboptimal performance.

The algorithm proceeds as follows.
The samples delivered to the laboratory are monitored from the perspective of their priority and arrival time.
The centrifuge is loaded with the samples and started afterward if the available samples fully meet the centrifuge capacity or if the designed timeout has expired after the last sample has arrived.
For this, two different timeouts are considered---a long one for statim/routine samples and a short one for vital samples.
The short timeout applies whenever a vital becomes available. 
The samples to be batched are chosen by the priority-first rule (i.e., all available vitals are loaded first, then statims, and finally routines), with ties broken by the sample arrival time (i.e., first in, first out).
We assume that the timeout is set to 2 minutes for vitals and 4 minutes for statim/routine samples~\cite{yang2015optimization}.

\begin{figure}[ht]
    \centering
    \begin{tikzpicture}
        \begin{axis}[
            height=0.25\textwidth,
            width=0.75\textwidth,
            ymin=0,
            ymax=0.5,
            xmin=0,
            xmax=25,
            xlabel={centrifuge start time $S^B$ [min]},
            label style={font=\scriptsize},
            tick label style={font=\scriptsize},
            minor x tick num=4,
            hide y axis
            ]
            \draw [->,>=stealth,blue] (axis cs:1,0.0) -- (axis cs:1,0.15);
            \draw [-,blue] (axis cs:1,0.15) -- (axis cs:8,0.15);
            \draw [->,>=stealth,blue] (axis cs:8,0.15) -- (axis cs:8,0.0);
            \node[rotate=0] at (axis cs:1,0.4){\scriptsize{$\hat{r}_j$}};
            \node[rotate=0] at (axis cs:8,0.4){\scriptsize{$\hat{a}_j$}};
    
            \draw [->,>=stealth,blue] (axis cs:3,0.0) -- (axis cs:3,0.3);
            \draw [-,blue] (axis cs:3,0.3) -- (axis cs:9,0.3);
            \draw [-,dashed,blue] (axis cs:9,0.3) -- (axis cs:9.5,0.3);
            \draw [->,>=stealth,dashed,blue] (axis cs:9.5,0.3) -- (axis cs:9.5,0.0);
            \node[rotate=0] at (axis cs:3,0.4){\scriptsize{$\hat{r}_{j+1}$}};
            \node[rotate=0] at (axis cs:9.5,0.4){\scriptsize{$\hat{a}_{j+1}$}};
            \draw[fill=gray!20] (axis cs:8,0.15) rectangle (axis cs:9,0.21);
            \draw [->,>=stealth,gray] (axis cs:10.5,0.25) -- (axis cs:9.1,0.18);
            \node[gray] at (axis cs:14.4,0.25) {\footnotesize{centrifuge timeout}};
        
            \draw[pattern=north west lines, pattern color=red] (axis cs:9,0.05) rectangle (axis cs:24,0.1);
            \node[red] at (axis cs:16.8,0.16){\footnotesize{centrifuging}};
            
            \draw [-,thick,red!40] (axis cs:9,0.00) -- (axis cs:9,0.35);
            \draw [-,thick,red!40] (axis cs:9,0.35) -- (axis cs:11,0.35);
            \node[color=red!40] at (axis cs:13.53,0.35){\footnotesize{current time}};
        \end{axis}
    \end{tikzpicture}
    \caption{Example demonstrating the unsuitable timing of the centrifuge. The centrifuge is started too early.}
    \label{fig:centrifuge_early}
\end{figure}
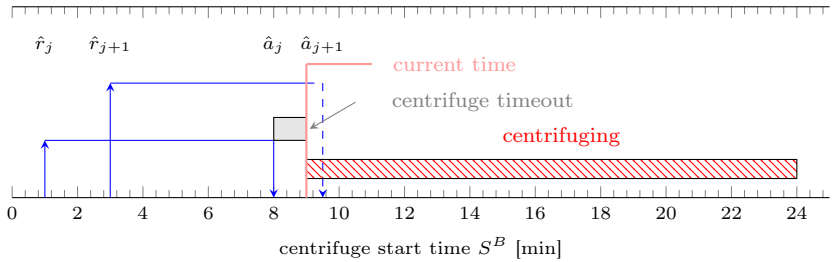
The main source of inefficiency of this strategy is the following.
The scenario in Figure~\ref{fig:centrifuge_early} represents the situation in which sample $j$ is available at time $8$ while sample $j+1$ is registered in the system but has not yet arrived.
If the centrifuge is activated after the timeout (e.g., $1$ minutes) past the arrival of sample $j$, then sample $j+1$ arrives just a few moments later and cannot be included in the current batch, and its TAT is prolonged by at least one spinning cycle.
Thus, even if the timeouts are designed to mitigate these situations, they cannot be completely avoided.
This shortcoming motivates us to design an algorithm that account for the samples in transport and is able to postpone the run of the centrifuge until the vital sample arrives.
This strategy is described in the following section.
\subsection{Look-ahead policy} \label{sec:fifo_alg}
To address the main disadvantage of the baseline solution (threshold-based policy), which is agnostic to the samples that are currently in transit, we propose a look-ahead batching algorithm.
Thus, in contrast to the algorithm described in Section~\ref{sec:sota_alg}, it attempts to prevent the situations illustrated in the example in Figure~\ref{fig:centrifuge_early}, where a vital sample is postponed because of an overly early centrifuge run.
This is achieved by utilizing information about samples registered in the system prior to their arrival at the laboratory.
All samples currently in transit are \textit{considered}, but the selection can be further refined by mitigating samples that were registered in the past $d$ minutes, for example.
The reason for this refinement is to steer the look-ahead policy toward samples that are likely to arrive soon, rather than blocking the centrifuge with samples that were registered just now and thus will not arrive anytime soon.

The rules governing the dispatch of the centrifuge and sample batching resemble natural requirements, prioritizing vital samples through the system.
Essentially, dispatch is guided by the following rules: (a)~if a vital sample is available and no other vital sample is currently under \textit{consideration}, then run the centrifuge, and (b)~if no samples are included in the centrifuge within the timeout and no vital sample is under consideration, then run the centrifuge, or (c) run the centrifuge if no samples are included in the centrifuge within the defined timeout and a vital sample is available.
The timeout is a parameter of the policy subject to modification.
In our experiments, we kept the same values as in the currently deployed threshold-based policy. 
The above rules aim to anticipate the registration of unforeseen vital samples to mitigate excessively early runs when centrifugation commences just before a vital sample arrives.

The look-ahead strategy is driven by both incoming and registered vital samples, with the goal of mitigating problems with the centrifuge starting too early and resulting in fewer vital sample misses.
However, even considering samples currently in transit without assuming their distribution of transportation times does not guarantee efficient control of the centrifuge.
The main difficulty with the system based on priority rules and hard-coded logic is that it is not flexible enough.
Given a fixed set of decision rules, one can often devise an adversarial scenario that these rules would not handle well.
The emergence of these corner cases is also amplified by the presence of many samples in the system with different distributions of transport times and processing time requirements, which leads to conflicting goals.
            
\begin{figure}[ht]
    \centering
        \begin{tikzpicture}
        \begin{axis}[
            height=0.26\textwidth,
            width=0.75\textwidth,
            ymin=0,
            ymax=0.9,
            xmin=0,
            xmax=25,
            xlabel={centrifuge start time $S^B$ [min]},
            label style={font=\scriptsize},
            tick label style={font=\scriptsize},
            minor x tick num=4,
            hide y axis
            ]
            \draw [->,>=stealth,blue] (axis cs:1,0.0) -- (axis cs:1,0.15);
            \draw [-,blue] (axis cs:1,0.15) -- (axis cs:4.5,0.15);
            \draw [->,>=stealth,blue] (axis cs:4.5,0.15) -- (axis cs:4.5,0.0);
            \node[rotate=0] at (axis cs:1,0.39){\scriptsize{$\hat{r}_j$}};
            \draw [->,>=stealth,blue] (axis cs:3.5,0.0) -- (axis cs:3.5,0.3);
            \node[rotate=0] at (axis cs:3.3,0.39){\scriptsize{$\hat{r}_{j+1}$}};
            \draw [-,blue] (axis cs:3.5,0.3) -- (axis cs:4.5,0.3);
            \draw [-,dashed,blue] (axis cs:4.7,0.3) -- (axis cs:18,0.3);
            \draw [->,>=stealth,dashed,blue] (axis cs:18,0.3) -- (axis cs:18,0.0);
            \node[rotate=0] at (axis cs:4.6,0.39){\scriptsize{$\hat{a}_j$}};
            \node[rotate=0] at (axis cs:18,0.39){\scriptsize{$\hat{a}_{j+1}$}};
            
            \addplot[smooth,fill=green,fill opacity=0.4]
            coordinates{(12,0) (16,0.03) (17,0.1) (18,0.15) (18.5,0.01) (20,0)};
            \addplot[smooth,fill=blue,fill opacity=0.3]
            coordinates{(10,0) (8.5,0.03) (8,0.1) (7,0.15) (6,0.01) (2,0)};
            \draw[pattern=north west lines, pattern color=red] (axis cs:18,0.05) rectangle (axis cs:29,0.1);
            \node[red] at (axis cs:21.5,0.16){\footnotesize{centrifuging}};
            \draw [-,thick,red!40] (axis cs:4.5,0.16) -- (axis cs:4.5,0.3);
            \draw [-,thick,red!40] (axis cs:4.5,0.62) -- (axis cs:4.5,0.5);
            \draw [-,thick,red!40] (axis cs:4.5,0.62) -- (axis cs:5,0.62);
            \node[color=red!40] at (axis cs:7.53,0.625){\footnotesize{current time}};
        \end{axis}
    \end{tikzpicture}
    \caption{Example demonstrating the unsuitable timing of a centrifuge with uncertain transport times. The centrifuge is started too late.}
    \label{fig:centrifuge_late}
\end{figure}
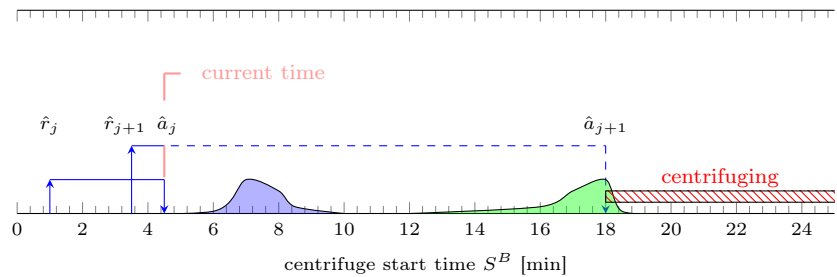

To see this, consider the following simple scenario in Figure~\ref{fig:centrifuge_late}. 
Sample $j$ is available at time $t=4.5$, while sample $j+1$ is in transit, similar to the scenario in Figure~\ref{fig:centrifuge_early}.
However, waiting for the arrival of sample $j+1$ could be a suboptimal decision in this case since the efficiency of the decision to wait for the transiting sample depends on the distribution of its transport times, as shown in green in Figure~\ref{fig:centrifuge_late}. In that case, it is more beneficial to run the centrifuge immediately and not wait for the arrival of sample $j+1$.
Sample $j+1$ is unlikely to arrive before the centrifuge finishes; thus, $j+1$ can be comfortably included in the next run, and the TAT of sample $j$ is unaffected.
However, if the transport distribution is like the one depicted in light blue, then it is better to postpone the centrifuge until $j+1$ arrives.

The examples in figures~\ref{fig:centrifuge_early} and \ref{fig:centrifuge_late} highlight some of the difficulties encountered when designing rule-based solutions that utilize the information about samples in transit but do not account for the probability distributions of the transport times.
Enhancing the set of decision rules with information about the distribution of transport times leads to an optimization problem rather than a fixed decision tree that triggers the centrifuge to run at a specific time on the basis of some designed conditions.
Therefore, a more flexible and informed solution is needed to solve this problem efficiently.

In the following sections, we propose to improve the timing of the centrifuge with a mathematical model that minimizes the total flow time of the available samples under a soft-deadline constraint with respect to the artificially chosen due date.
To include the uncertainty in the transport times, we also minimize the expected flow time of samples in transit.
Combining these two terms inside the objective function, the model provides centrifuge timings requested by the available samples while also accounting for the arriving samples according to the estimated transport time distributions.

The uncertainty in our optimization problem is imposed via uncertain registration times $\hat{r}_j$, transport times $\hat{\tau}_j$, priorities, and the total number of samples.
In the considered setting at work, we do not assume any particular distribution over registration times $\hat{r}_j$; therefore, we treat them as unpredictable random events.
Thus, the only uncertain parameters reflected in our model are the collection of transport times $\{\hat{\tau}_j\}_{j}$ of samples currently in transport.

General strategies for solving optimization problems with uncertain parameters can be classified into two main groups.
The first group of methods aims to replace the uncertain parameter with its estimate, hence reducing it back to a deterministic problem.
The most typical choice might be to use expected values, some quantile of the empirical distribution, or essentially any other single-point estimate. 
This intuitive approach works well if the underlying distribution of the uncertain parameter can be expressed with a single number sufficiently well, e.g., for a normal distribution with low variance.

The second group of methods models the uncertainty of the parameters as a part of the optimization model---typically expressed via an uncertainty set or by a cumulative distribution function (CDF).
These methods are especially useful in cases in which the distribution cannot be effectively expressed via single-point estimates, e.g., complex distributions with high variance or even multimodal distributions.
If the transport time distribution took a bimodal form with two peaks, as shown, e.g., in Figure~\ref{fig:centrifuge_late} in blue and green, then no single-point estimate would effectively capture the nature of such a distribution.
On the other hand, if the distribution of the parameter is narrow enough, then this stochastic modeling essentially converges to a single-point method and thus may be viewed as a more general method.

In the considered problem with sample transportation, it is reasonable to assume that the resulting distribution could be complex, e.g., because it reflects time-of-day factors or even because of the possibility of different modes of transportation (courier or tube mail).
Therefore, in the following section, we propose a stochastic batching mathematical model that models the transport time uncertainty via a CDF as a part of the optimization model.

\subsection{Stochastic batching mathematical model}\label{sec:mip_alg}
In Section~\ref{sec:fifo_alg}, we have demonstrated how the unsuitable timing of the centrifuge can negatively impact the values of patient TAT.
To mitigate timings that are too early and too late, we propose a stochastic batching mathematical model that exploits information about the registration of new samples in the system and the estimated transport time distribution.
The mathematical model is included in the DES simulation and is used to determine the timing of the next centrifuge run and the content of the batch.

The core principle of the model is to prioritize vital samples exclusively, i.e., minimize the overall patient TAT of vital samples. Consequently, the centrifuge's operational schedule is determined by the processing requirements of these vital samples. Nevertheless, the capacity of the centrifuge is shared with statim and routine samples, which are loaded in a priority-based first-in, first-out order.
As demonstrated in Section~\ref{sec:experiments}, this approach enables real-time operation of the model without adversely affecting the TAT of non-vital samples.

We express the criterion as the minimization of the sum of the expected flow times, i.e., the expected completion times of samples minus their release times~\cite{NOVAK2022438}, with a soft deadline constraint.
The reason we do not directly optimize, e.g., for the 0.95 quantile of the TAT distribution, is that it is not a criterion but rather a statistical indicator to keep track of, as it indicates whether an excessive number of samples are late in the presence of outliers introduced by random errors in the process.
Since we cannot affect the sources of these errors, the model optimizes the expected patient TAT.

The batching model performs relatively short-term decision-making because of the typical transport times of vital samples and centrifuge capacity. 
Specifically, we consider only which vital samples currently available or transported should be included in the next centrifuge run and what can be postponed until the following run.
The model distinguishes two sets of samples.
First, a set of vital samples that arrived at the laboratory and thus are available at time $t$ is denoted as $V_t$.
Next, it keeps track of the samples $\overline{V_t}$ that are registered in the system at time $t$ but have not yet arrived in the laboratory.
For the available samples $j\in V_t$, we possess the full information; thus, we can express the time $C_j$ when sample $j$ is completed as a function of the centrifugation start time.
As suggested above, one of the main functions of the model is to determine the content of the batches for the next two consecutive runs of the centrifuge.
This is based on the assumption of typical transportation times for vital samples in our data, the centrifuge's capacity, and its cycle time.
With those, every transiting vital sample can be processed within two upcoming spinning cycles.
Thus, the model uses optimization variables $S^B_1$ and $S^B_2$ to determine the start times of the next batch or the second next batch, respectively.
Thus, we can compute $C_j$ depending on whether the sample is assigned to the current or the next batch using binary variables $y_{j,1}$ and $y_{j,2}$.
An available vital sample might be postponed to the next batch (i.e., $y_{j,2}=1$) because of the limited capacity $N$ of the centrifuge, although in practice, this would unlikely be a binding constraint.

On the other hand, we cannot compute the completion time of samples that are in transit because of the uncertainty of the transport time.
Therefore, for samples $j\in\overline{V_t}$, we compute the probability that sample $j$ will arrive before the scheduled start of the $k$-th centrifuge batch, denoted as $\pr{j}{k}$ with $k\in\{1,2\}$, as the considered horizon assumes two upcoming batches.
This probability is computed from the cumulative distribution function of the transport time of sample~$j$, which is estimated from historical data as a function of the sample priority $\lambda_j$ and the ward $w_j$. The use of an empirical distribution may capture potentially complex transport time distributions, enabling efficient centrifuge scheduling by identifying low-probability arrival windows.
Modeling transport time uncertainty via its full CDF may be a more effective choice, for example, when a bimodal distribution for which a single-point estimate, such as the expected value, could produce an outcome that even has a zero probability of occurrence.
%

A natural question arises regarding whether the distribution modeled via its CDF is stationary; that is, whether it depends on the particular hour of the day. The model we propose below does not make any assumptions about this.
Therefore, the distribution can even be nonstationary and can be incorporated by simply using a different CDF for each hour of the day if enough data are available.
In this work, we assume that the distribution does not change over the course of a day; however, the model may benefit from utilizing more of the structure hidden within the transport time distribution.

\begin{figure}[ht]
    \centering
    \begin{tikzpicture}
        \begin{axis}[
            height=0.36\textwidth,
            width=0.8\columnwidth,
            xlabel={$S^B_1-\hat{r}_j$ [min]},
            ylabel=$\pr{j}{1}$,
            label style={font=\scriptsize},
            tick label style={font=\scriptsize},
            ]
            \addplot[smooth,thick,dashed] table {cdf_ortjip.txt};
            \addplot[thick,blue,draw opacity=0.6] table {cdf_ortjip_crop.txt};
            \addplot[thick,blue,dashed,draw opacity=0.6] table {cdf_ortjip_crop_start.txt};
            
            \draw [dashed] (axis cs:11,1) -- (axis cs:11,-0.1) node[above=6.2em,rotate=90] {\footnotesize{$j$ cannot be in the batch}};
            \draw [->,>=stealth] (axis cs:11,1) -- (axis cs:0,1);
            \draw [dashed] (axis cs:64,1) -- (axis cs:64,-0.1) node[below=-6em,rotate=90] {\footnotesize{$j$ can be in the batch}};
            \draw [->,>=stealth] (axis cs:64,0.03) -- (axis cs:75,0.03);
            \draw[->,>=stealth,blue] (axis cs:35,0.15) node[right] {\scriptsize{$t-\hat{r}_j$}} -- (axis cs:28.5,0.28) ;
        \end{axis}
    \end{tikzpicture}
    \caption{Probability that sample $j$ from \textsc{ortjip} department registered at time $\hat{r}_j$ can be included in the upcoming batch as a function of batch start time $S^B_1$.}
    \label{fig:cdf-example}
\end{figure}
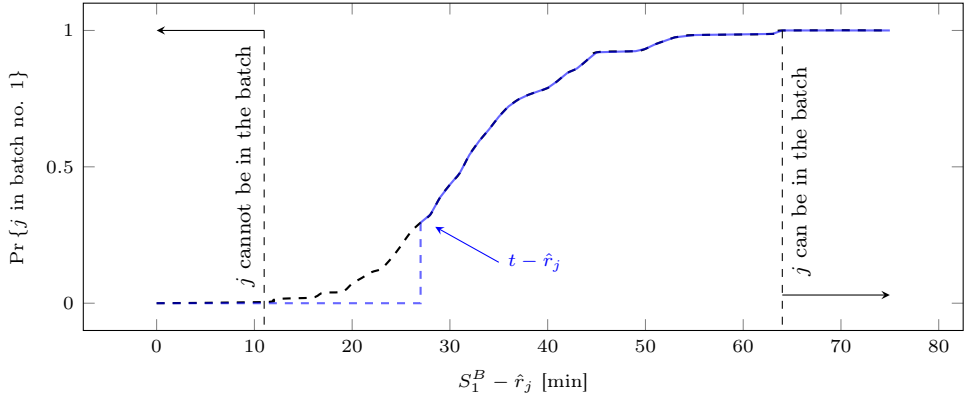

The example in Figure~\ref{fig:cdf-example} shows how the probability of sample $j$ being included in the next centrifuge run depends on the centrifuge start time $S^B_1$.
If~$S^B_1$ is too small, then sample $j$ will surely miss the run because of its minimum transport time.
When it is large enough, then sample $j$ will surely be available in the laboratory to be centrifuged.
As simulation time $t$ progresses, the relationship between the centrifuge start time and the probability of sample arrival is readjusted such that $S^B_1\geq t$ always holds, as shown in blue in Figure~\ref{fig:cdf-example}.

Although we compute the probability that a sample arrives before the start of the centrifuge for each vital sample currently in transport, this does not necessarily mean that the model decides to include all these samples in the earliest possible batch.
Similarly, as in the case of available samples $V_t$, the model can postpone transiting vital sample $j$ to the next batch.
Thus, for vital samples $j\in\overline{V}_t$ that are currently in transport, we introduce a binary variable $\overline{y}_{j,2}$ with meaning $\overline{y}_{j,2}=1$ when transiting sample $j\in\overline{V}_t$ is considered for the second batch.
Allowing this degree of freedom is beneficial, for example, in cases where the minimum transport time of the transiting sample is greater than the centrifuge cycle time~$\delta$.

To model this decision, we introduce a new set of variables $\prp{j}{k}$, which models the probability that sample $j\in\overline{V_t}$ will be included in the batch $k\in\{1,2\}$. 
The difference from the $\pr{j}{k}$ variables is that it reflects the value of $\overline{y}_{j,2}$; thus, if sample $j\in\overline{V_t}$ is set to be postponed to the second centrifuge run (i.e., $\overline{y}_{j,2}=1$), then $\prp{j}{1}=0$.
These variables are used to express a part of the criterion that computes the expected flow time of samples in transport. 
The full model is given by equations \eqref{eq:model_obj}--\eqref{eq:model_end}.
The objective \eqref{eq:model_obj} minimizes the sum of the flow times for the available samples plus the expected flow time of the samples in transit.
The last term of the objective penalizes the use of $\xi_j$ and $\overline{\xi}_j$ variables for violation of the soft deadline constraint.
Constraints \eqref{eq:model_sb1_constr} and \eqref{eq:model_sb2_constr} model the minimum possible start times of the next two following runs of the centrifuge.
The constraint~\eqref{eq:model_sb1_constr} asserts that the next start time must be greater than or equal to the last centrifuge start time $S^B_0$ plus its cycle time $\delta$.
The constraint~\eqref{eq:model_sb2_constr} defines the ordering of the first and the second considered batches.
The expected utilization of centrifuge capacity is constrained by \eqref{eq:model_capacity1}.
The constraint~\eqref{eq:model_batch_selection} enforces that each available sample is included in either the first or the second centrifuge run, while constraints \eqref{eq:model_completion1} and \eqref{eq:model_completion2} are indicator constraints used to compute the completion time $C_j$ of sample $j\in V_t$ based on the start time of the assigned batch and the processing time~$p_j$ of the sample in the analytic part of the system.

The constraint \eqref{eq:model_cdf} models the probability of an early sample arriving as a piecewise linear approximation of the sample transportation time cumulative distribution function $\text{CDF}_{j}(x) \in [0,1]$, exactly as demonstrated in Figure~\ref{fig:cdf-example}.
The assumption that a vital sample currently in transit will surely arrive for the second batch at the latest is expressed by \eqref{eq:model_pr_sum}. 
The indicator constraints \eqref{eq:model_prp_start}--\eqref{eq:model_prp_end} interconnect the semantic of \textit{the probability $\pr{j}{k}$ that sample $j$ can be included} and the actual planned decision $\overline{y}_{j,2}$ that \textit{sample $j$ is included} in the corresponding batch.
%
Finally, the constraints \eqref{eq:model_soft_deadline_start}--\eqref{eq:model_soft_deadline_end} impose a soft deadline in that the (expected) completion times of both available and transiting vital samples must be before the requested due date of $D$ seconds from registration via $\xi_j, \overline{\xi}_j$, slack variables.

We include this constraint since, in preliminary experiments, we observed that the model occasionally yields solutions that sacrifice the worst-case performance (i.e., maximum patient TAT) in favor of optimizing the lower quantiles.
However, including this constraint is solely the user's choice and may be omitted if a small percentage of high TAT values does not represent an issue.

The model resembles a \mbox{nonconvex} \mbox{mixed-integer} quadratic program (MIQP) but can be solved up to global optimality by modern solvers, such as Gurobi.
The full model is given as follows:
\begin{align}
    &\min \sum_{j\in V_t}  \left(C_j - \hat{r}_j\right)+\sum_{k\in\{1,2\}}\sum_{j\in\overline{V_t}}\prp{j}{k}\cdot\left(S^B_k+\delta + p_j-\hat{r}_j\right)+\notag\\
    &+M\cdot \left(\sum_{j\in V_t} \xi_j + \sum_{j\in\overline{V}_t} \overline{\xi}_j \right)\label{eq:model_obj}\\
    \intertext{subject to}
    &S^B_1 \geq \max\{t, S^B_0 + \delta\} \label{eq:model_sb1_constr}\\
    &S^B_2 \geq S^B_1 + \delta\label{eq:model_sb2_constr}\\
    &\sum_{j\in V}y_{j,1}+\sum_{j\in\overline{V_t}}\prp{j}{1}\leq N\label{eq:model_capacity1}\\
    &y_{j,1}+y_{j,2} = 1\quad \forall j \in V_t\label{eq:model_batch_selection}\\
    &y_{j,1} \implies C_j = S^B_1 + \delta + p_j\quad \forall j \in V_t\label{eq:model_completion1}\\
    &y_{j,2} \implies C_j = S^B_2 + \delta + p_j\quad \forall j\in V_t\label{eq:model_completion2}\\
    &\pr{j}{1}= \text{CDF}_{j}(S^B_1-r_j) \quad \forall j\in \overline{V_t}\label{eq:model_cdf}\\ 
    &\pr{j}{1}=1-\pr{j}{2}\quad \forall j\in\overline{V_t} \label{eq:model_pr_sum}\\
    &\overline{y}_{j,2} = 1 \implies \prp{j}{1} = 0\quad \forall j \in \overline{V_t}\label{eq:model_prp_start}\\
    &\overline{y}_{j,2} = 1 \implies \prp{j}{2} = 1\quad \forall j \in \overline{V_t}\\
    &\overline{y}_{j,2} = 0 \implies  \prp{j}{1} =\pr{j}{1}\quad \forall j \in \overline{V_t}\\
    &\overline{y}_{j,2} = 0 \implies \prp{j}{2} =\pr{j}{2}\quad \forall j \in \overline{V_t}\label{eq:model_prp_end}\\
    &C_j - \hat{r}_j \leq D + \xi_j \quad \forall j \in V_t \label{eq:model_soft_deadline_start}\\
    &\sum_{k\in\{1,2\}}\sum_{j\in\overline{V_t}}\prp{j}{k}\cdot(S^B_k+\delta + p_j-\hat{r}_j) \leq D + \overline{\xi}_j \quad \forall j \in \overline{V}_t\label{eq:model_soft_deadline_end}\\
    \intertext{where}
    &S^B_0 = \text{the last centrifuge start time}\\
    &y_{j,1},y_{j,2}\in \{0,1\}\quad \forall j\in V_t\\
    &\overline{y}_{j,2}\in \{0,1\}\quad\forall j \in \overline{V_t}\\
    &S^B_1,S^B_2\geq 0\\
    &C_j\geq0 \quad \forall j\in V_t\\
    &\pr{j}{k}\in [0,1] \quad \forall j \in \overline{V_t}, k \in \{1,2\}\\
    &\prp{j}{k}\in [0,1] \quad \forall j \in \overline{V_t}, k \in \{1,2\}\\
    &\xi_j \geq 0 \quad \forall j\in V_t\\
    &\overline{\xi}_j \geq 0 \quad \forall j\in\overline{V}_t \label{eq:model_end}
\end{align}

The model \eqref{eq:model_obj}--\eqref{eq:model_end} is employed inside the online setting in the following way.
As asserted above, the stochastic MIQP model is driven only by vital samples.
Although the policy is agnostic to the registration and arrival of statim and routine samples, it still shares the available capacity of the centrifuge with other available samples.
At time instances $t$ of the DES simulation when a vital sample is present in the system (i.e., $V_t\neq\emptyset$), the model \eqref{eq:model_obj}--\eqref{eq:model_end} is constructed and solved.
Suppose that the (optimal) value of the $S^B_1$ variable is greater than $t$. In that case, the model registers a custom wake-up event (see Figure~\ref{fig:DES}, Policy process) for a future time instant equal to the computed value of the $S^B_1$ variable.
If the value of the $S^B_1$ variable equals $t$, then the centrifuge is dispatched with the vital samples $j\in V_t$, for which $y_{j,1}=1$.
If this amount is less than the centrifuge capacity $N$, then the rest of the capacity is filled with available samples in order of their priority, with ties broken by the sample arrival time.
When no vital sample is available or registered in the system, the model \eqref{eq:model_obj}--\eqref{eq:model_end} is not constructed. 
Instead, the centrifuge is operated following the look-ahead policy described in Section~\ref{sec:fifo_alg}.

Here, we would like to make a few remarks about the model \eqref{eq:model_obj}--\eqref{eq:model_end}.
First, note that the objective function \eqref{eq:model_obj} is the sum of (expected) flow times plus the soft constraint violation penalty.
This choice seems to be natural, as the flow time $F_j = C_j - \hat{r}_j$ is precisely the patient TAT of a sample $j$.
However, we also experimented with different objective functions, such as $\sum T_j$ with $T_j=\max\{0, C_j - d_j\}$, where $d_j$ is a virtual due date set as $d_j=\hat{r}_j + \beta$ for some $\beta > 0$.
The idea is to allow each sample to accommodate some defined lateness $\beta$ (e.g., 45 minutes) before incurring a penalty.
With that, the model becomes more flexible in postponing the centrifuge run to accommodate more transiting samples.
However, with our data, the total tardiness or other related objectives, such as minimizing the number of tardy jobs $\sum U_j$ or minimizing squared flow times, did not present better performance than the total flow time minimization, as will be shown in Section~\ref{sec:exp_objectives}.

Finally, we highlight the limitations of the above model, which utilizes a full CDF transport time distribution.
Specifically, one may wonder whether modeling a full CDF of the transport time distribution yields a significant improvement over simpler (e.g., look-ahead) methods, which may be equipped with just single-point estimates of transport times. The first such setting is where the transport time distribution is narrow or, at the very least, close to a deterministic distribution that can be efficiently captured by a single parameter. 
This situation is discussed at the end of Section~\ref{sec:fifo_alg}.
Other scenarios where the knowledge of the full CDF does not present a significant advantage include cases where either too few or too many vital samples are released during the day.
In the first scenario, the method effectively converges toward a look-ahead policy, whereas in the second scenario, the vital priority essentially becomes statim, as there are not enough lower-priority samples to take advantage of.
Additionally, if vital samples tend to have a long transport time relative to the centrifuge cycle time, then the relative improvement in patient TAT tends to diminish.
We study this case in Section~\ref{sec:system_params}, where we investigate the relative improvement over the look-ahead method for scenarios with different cycle times~$\delta$.
\subsection{Offline perfect-knowledge algorithm}%
\label{sec:perfect_knowledge}%
To construct a point of reference for the proposed online algorithms, we also consider a deterministic, offline variant of the problem, with full information about uncertain parameters available to the decision-maker. That is, for the offline algorithm, the release times of the samples and the actual values of the transportation times are known in advance. In this sense, the results of such a perfect knowledge algorithm represent the theoretical boundary that can be attained by any method provided that the offline algorithm is exact. 
First, we describe a general schema of the proposed offline algorithm and discuss the design choices. 
Then, we provide details of two algorithms for different criteria functions that can be incorporated into the introduced optimization schema.

Baptiste et al.~\cite{Baptiste2000} reported that $1\ |\ p\text{-batch},\ b<n,\ r_j,\ p_j=p\ |\ f$ is polynomial for at least $f \in \{ \sum w_j U_j,\ \sum w_j C_j,\ \sum T_j,\ \max_j T_j\}$ and some other objective functions following specific properties. 
While these results could also be applied to solve some variants of our problem in an offline perfect-knowledge setting, their algorithm has time complexity in $O(n^8)$. 
However, a typical laboratory processes hundreds to thousands of samples daily, which renders any high-order polynomial algorithm intractable in practice. For smaller instances, e.g., when only vital samples are considered, polynomial algorithms can indeed be used, but such small instances can also be easily solved by modern MIP solvers.

While minimizing functions such as $\sum w_j U_j,\ \sum w_j C_j,\ \sum T_j$ or $\max_j T_j$ for instances beyond the reach of MIP solvers may be relatively difficult because of the high polynomial order of the algorithm of \cite{Baptiste2000}, the feasibility problem $1\ |\ p\text{-batch},\ b < n,\ r_j,\ p_j = p,\ \text{deadlines } d^\prime_j\ |-$ with deadlines $d^\prime_j$ is solvable in $O(n^2)$ time when the algorithm proposed in~\cite{Condotta2010} is used.
There, an $O(n^2 \log n)$-time algorithm was proposed for $1\ |\ p\text{-batch},\ b < n,\ r_j,\ p_j = p\ |\ L_{\text{max}}$, i.e.\ minimizing the maximum lateness $\max_j L_j = \max_j C_j - d_j$ for the case when $d_j$ are soft bounds (i.e., due dates).
A similar approach can be used to minimize the maximal patient TAT while respecting the priority of samples, as will be shown later.

We designed an offline perfect-knowledge algorithm to address the challenges described above, taking advantage of the properties of the problem; that is, the hierarchical properties of sample priorities---a vital sample is more important than any statim sample, and the same relation holds for statim and routine samples. 
Therefore, we can use a hierarchical approach, where we first derive completion times for vital samples alone without considering other sample priorities. 
Then, the algorithm schedules statim samples and, finally, the routine ones so that the completion times of vitals and statims derived by the previous stages are not violated. 
We refer to the problem of scheduling samples of a certain priority as a stage. 
Each stage may use a different objective function and is solved by its own algorithm, which we refer to as a \textit{subroutine}.

Details are described in Appendix~\ref{app:offline}. 
The pseudocode of the hierarchical optimization method is shown in Algorithm~\ref{alg:revised_full_algorithm}. 
The algorithm starts with the subroutine \texttt{SubroutineVital}. Each subroutine accepts three parameters: a set of samples $\mathcal{P}$ scheduled in the previous stage, whose completion times are $C^P$, where $C^P_i$ is the completion time of sample $i \in \mathcal{P}$, and a set of samples to be scheduled during the current stage $\mathcal{J}, \ \mathcal{P} \subseteq \mathcal{J} \subseteq \mathcal{S}$.
The subroutine provides one or more (optimal) schedules for the samples contained in $\mathcal{J}$, each characterized by completion times. 
Here, in line 2, the set of schedules returned by the subroutine is denoted by $\mathcal{C}$.
At the beginning, the set $\mathcal{J}$ consists exclusively of vital samples.
In the first stage, there are no previously scheduled samples $\mathcal{P}$; thus, their completion times are not available. 
Then, statim samples are scheduled by \texttt{SubroutineStatim}, with $\mathcal{J}$ consisting of both vital and statim samples. 
Since it is the second stage, $\mathcal{P}$ contains all the vital samples, and each schedule from $\mathcal{C}$ must be considered separately (lines 4--5).
Therefore, the subroutine is executed for each schedule obtained from the \textsc{vital} stage, and the results are combined into a single set of new schedules $\mathcal{C}'$. 
Next, all the schedules with the optimal value of the objective function from the \textsc{statim} stage are preserved and passed into $\mathcal{C}$. The same steps are followed for \textsc{the routine} samples. 
Finally, the start times of the batches $S^B$ are calculated on the basis of one of the optimal schedules from $\mathcal{C}$. 
If the subroutines at each stage are exact and return all optimal schedules, the hierarchical algorithm is also exact. 
We introduced two subroutines, built upon exact algorithms: the total TAT minimization subroutine (based on an MIP formulation) and the maximum TAT minimization subroutine (based on a polynomial-time algorithm). 
Both are described in detail in Appendix~\ref{app:offline}. 
\begin{algorithm}[ht]
\DontPrintSemicolon
\caption{Offline perfect-knowledge algorithm.}%
\label{alg:revised_full_algorithm}
\SetKwFunction{Subroutinea}{SubroutineVital}
\SetKwFunction{Subroutineb}{SubroutineStatim}
\SetKwFunction{Subroutinec}{SubroutineRoutine}
\SetKwProg{Fn}{Function}{:}{\KwRet}
\KwData{Samples $\mathcal{S}$.}
\KwResult{Start times of batches $S^B$.}
\tcp{Vital stage}
$\mathcal{P}\leftarrow \emptyset,\ \mathcal{J} \gets \{ j \in  \mathcal{S}\ |\ \lambda_j = \textsc{vital} \}$\;
$\mathcal{C} \leftarrow \Subroutinea{\ensuremath{\emptyset, \mathcal{P}, \mathcal{J}}}$\;
\tcp{Statim stage}
$\mathcal{C}' \leftarrow \emptyset,\ \mathcal{P}\leftarrow \mathcal{J},\ \mathcal{J} \gets \{ j \in  \mathcal{S}\ |\ \lambda_j = \textsc{statim} \} \cup \mathcal{J}$\;
\ForEach{$C^P \in \mathcal{C}$}{
    $\mathcal{C}' \!\!\leftarrow \mathcal{C}'\! \cup \Subroutineb{\ensuremath{C^P\!\!, \mathcal{P}\!, \mathcal{J}}}$\;
}
$\mathcal{C} \leftarrow \text{select optimal schedules from } \mathcal{C}'$\;
\tcp{Routine stage}
$\mathcal{C}' \leftarrow \emptyset,\ \mathcal{P}\leftarrow \mathcal{J},\ \mathcal{J} \gets \{ j \in  \mathcal{S}\ |\ \lambda_j = \textsc{routine} \}  \cup \mathcal{J}$\;
\ForEach{$C^P \in \mathcal{C}$}{
    $\mathcal{C}' \!\!\leftarrow \mathcal{C}'\! \cup \Subroutinec{\ensuremath{C^P\!\!, \mathcal{P}\!, \mathcal{J}}}$\;
}
$C^P \leftarrow \text{one of the optimal schedules in } \mathcal{C}'$\;
$S^B \gets \textrm{sort}(\{C^P_j - \delta - p_j\ | \ j \in \mathcal{S}\})$\;
\end{algorithm}
\section{Experiments} \label{sec:experiments}
All the algorithms in the experiments were implemented in Julia~1.11.5.
The discrete event system (DES) simulation uses the \emph{ConcurrentSim} library, a Julia implementation of DES inspired by the well-known Python DES framework \emph{Simpy}.
For the solutions of the MIP and MIQP models, we use the Gurobi solver version 11.
The benchmarks were conducted on commodity hardware, i.e., an AMD Ryzen 9 5950X 16-core processor with 64 GB of RAM.
The source codes are available on GitHub\footnote{\url{https://github.com/CTU-IIG/centrifuge-batching}}.

In Section~\ref{sec:data_desc}, we describe which data were used for this case study and what data cleansing steps were applied.
We subsequently compare the distributions of patient TATs obtained with the methods discussed in this paper, utilizing different environmental information, as denoted in Table~\ref{tab:alg_list}. The aim of this experiment is to show the benefit of the proposed solution with respect to the baseline solution that is often used in laboratories.
Different objective functions and their effects on the solution are analyzed in Section~\ref{sec:exp_objectives}.
Next, we use the perfect-knowledge offline method introduced in Section~\ref{sec:perfect_knowledge} to establish the best possible results attainable with the current laboratory hardware.
Furthermore, we use our DES simulation and batching algorithms to explore the space of centrifuge parameters to assess under which configurations the proposed solutions achieve the best performance over the look-ahead policy.
These considerations are reflected in the managerial insights in Section~\ref{sec:insights}, which discusses the impact on the quality of the provided health care and the economic implications depending on the system configuration and the dispatching algorithm used.
\subsection{Data description} \label{sec:data_desc}
Our data for the purposes of this study come from University Hospital Královské Vinohrady, Czech Republic.
We obtained anonymized patient records for October 2018 from the Laboratory Information System (LIS), which were derived from the laboratory workflow analysis we performed for this hospital.

This work focuses on improving laboratory systems with a high level of automation, such as an automated centrifuge, sample distribution system, and analytical section, e.g., DxA 5000. Nevertheless, we do not concentrate on any specific laboratory system in the experimental analysis since the described approach is general and can be applied to other automated systems as well. On the other hand, since the data we use in our experiments are associated with the laboratory analyzers Cobas c702, c501, and e602, the processing time of each laboratory method, e.g., glucose, albumin, and troponin-I, assigned to each sample is set according to the parameters of those analyzers. Therefore, the processing time is deterministic and takes $p_j\in \{9, 10, 18, 27\}$ minutes depending on the prescribed method.
If the sample contains more methods to be performed, we estimate the processing time of the sample as the longest method it contains (methods are typically evaluated in parallel inside the laboratory analyzer). In total, 8799 routine, 6635 statim, and 142 vital samples were initially considered.

We excluded samples whose transportation time exceeded 200 minutes or whose duration between sample admission and result validation exceeded 24 hours. The procedure left 4936 routine, 4787 statim, and 137 vital samples.
Furthermore, for the vital samples, the transport time threshold was set to 22 minutes. Any transport time above this threshold was redrawn from $U(5, 22)$ to keep the overall number of vital samples intact while providing a more balanced input for the data generation procedure described below. These thresholds were selected on the basis of laboratory guidance to eliminate erroneous data, special samples processed only on specific days, and other outliers. Next, to ensure sufficient data for estimating transport time distributions in subsequent steps and to exclude insignificant, dummy, or outlier wards, we included only samples from the 30 hospital wards with the highest volume of nonvital samples while retaining all vital samples.

To enrich the dataset, we created alternative realizations of the uncertain parameters observed in the real data, resembling different variants of the observed days. That is, we generate new data by drawing realizations of transport times $\hat{\tau}_j$ from the estimated transport distribution of samples from the hospital ward $w$ associated with sample~$j$. If the number of samples from ward $w$ was insufficient for reliable estimation (fewer than 3 samples), we constructed a distribution using the combined sample set from such wards. This occurred only for vital samples, since otherwise, we kept only data from the top 30 wards.
In the case of routine and statim samples, the transport times are sampled from kernel density estimates derived from real data.
For the vital samples, owing to the small amount of data when conditioning on the hospital ward, the transport time distribution for hospital ward $w$ was estimated using a uniform distribution $\tilde{\tau}_j\sim U(lb_w, ub_w)$, where $lb_w$ and $ub_w$ are the lowest and the highest observed transport time realizations, respectively, for ward $w$ in real data.
This choice leads to a distribution with the maximum discrete entropy, thus imposing the least amount of bias and a priori knowledge on the system.
For the experiments below, we use the enriched dataset, which is worth 100 months of data.

\subsection{Patient TAT minimization}\label{sec:exp_patient_tat}
In this section, we compare the efficiency of the studied and proposed methods for centrifuge batching (summarized in~Table~\ref{tab:alg_list}).
Specifically, we compare the fixed-schedule policy, which serves as a deterministic and predictable baseline, with the threshold-based policy, which represents the solution often deployed in hospital laboratories~\cite{yang2015optimization}.
Next, we analyze the performance of the online methods proposed in this paper, which are the look-ahead policy from Section~\ref{sec:fifo_alg} and the stochastic MIQP model introduced in Section~\ref{sec:mip_alg}. All methods are executed in the online setting and interact with the environment using the DES simulation described in Section~\ref{sec:des}.

\begin{figure}[ht]
    \centering
    \pgfplotsset{
        cycle list/PuBu-8,
        cycle multiindex* list={
            mark list*\nextlist
            PuBu-8\nextlist
        },
    }

    \begin{tikzpicture}
        \begin{axis}[
            height=0.3\textwidth,
            width=0.6\textwidth,
            legend pos={outer north east},
            legend cell align=left,
            xmajorgrids, ymajorgrids,
            xlabel={Patient TAT [min]},
            ylabel={vital CDF [-]},
            label style={font=\scriptsize},
            tick label style={font=\scriptsize},
            table/x=probability,
            cycle list={blue!80, red!80, orange!80, cyan!80}
            ]
            \foreach \y in {miqp_cond_hiear_dd0,sota_vital_priority_with_timeout_algorithm,sota_algorithm,run_afap_algorithm}{
                \addplot+[thick,mark=none] table [y=\y, col sep=comma]{s52_ecdf.txt};
            }
            \addlegendentry{stochastic MIQP}\addlegendentry{look-ahead policy}\addlegendentry{threshold-based policy}\addlegendentry{fixed-schedule policy}
        \end{axis}
    \end{tikzpicture}
    \caption{Distribution of patient TAT for vital samples.}\label{fig:exp_16_comp_vitals}
\end{figure}
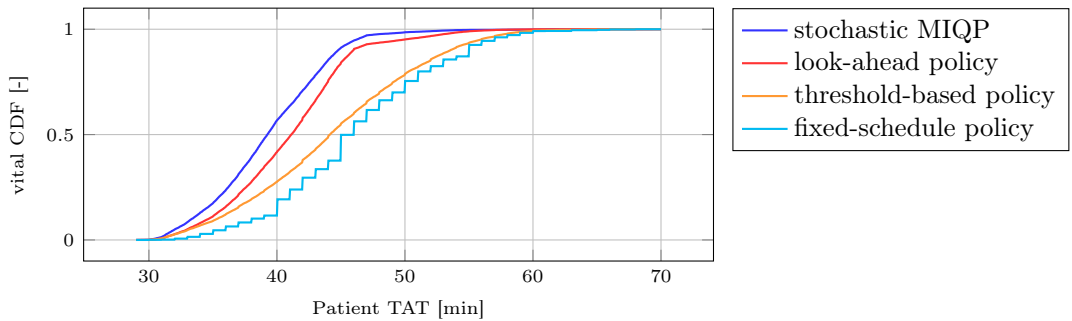

\begin{table}[ht]
    \centering
     \caption{Numerical results of patient TAT and laboratory TAT achievable by different policies.}\label{tab:ex_numerical_all}
    \begin{tabular}{*{7}{l}}
        \toprule
        Priority & Algorithm &  \multicolumn{2}{l}{Patient TAT [min]} & \multicolumn{2}{l}{Laboratory TAT [min]} \\
        \cmidrule(lr){3-4}\cmidrule(lr){5-6}
        & & 0.95 quantile & 0.5 quantile  & 0.95 quantile & 0.5 quantile \\
        \midrule
        vital    & stochastic MIQP   &   \textbf{46.2}  & \textbf{39.2} & \textbf{33.0} & \textbf{25.0}  \\
                 & look-ahead        &   49.9           & 41.1          & 36.5          & \textbf{25.0}  \\
                 & threshold-based   &   55.9           & 44.1          & 39.2          & 30.5  \\
                 & fixed-schedule    &   57.0           & 46.0          & 39.5          & 32.8  \\[.5em]
        statim   & stochastic MIQP   &  113.4           & 71.7          & 43.2          & \textbf{33.0}  \\
                 & look-ahead        &   \textbf{113.2} & 71.5          & 42.4          & \textbf{33.0}  \\
                 & threshold-based   &   \textbf{113.2} & 71.4          & \textbf{42.0} & \textbf{33.0}  \\
                 & fixed-schedule    &   115.0          & \textbf{70.0} & 43.0          & 33.7  \\[.5em]
        routine  & stochastic MIQP   &   129.3          & 87.6          & 49.6          & 35.5  \\
                 & look-ahead        &   129.1          & 87.4          & 48.4          & \textbf{35.4}  \\
                 & threshold-based   &   \textbf{129.0} & 87.2          & \textbf{48.0} & \textbf{35.4}  \\
                 & fixed-schedule    &   130.0          & \textbf{87.0} & \textbf{48.0} & \textbf{35.4}  \\
        \bottomrule
    \end{tabular}
   
\end{table}

In the first experiment, we compare the patient TAT distributions for all the methods.
The results are shown in Figure~\ref{fig:exp_16_comp_vitals}.
Generally, the more information an algorithm uses, the better the results that can be obtained for vital samples.
An important finding is that the stochastic MIQP and look-ahead policies are driven only by vital samples (i.e., their decision to trigger the centrifuge is not influenced by the arrival or registration of a nonvital sample) and do not incur significant performance degradation for those priorities, as shown in Table~\ref{tab:ex_numerical_all}.
The detailed results in the table reveal that stochastic MIQP achieves 3.7 minutes of improvement over the rule-based look-ahead policy in terms of the 0.95 quantile for vital samples, with the price of only 0.2 minutes degradation for statim samples.
With respect to the median performance, the patient TAT of vital samples is improved by 4.9 minutes when stochastic MIQP is used instead of the threshold-based policy.
In both cases, the improvement for vital samples was carried out only at minor degradation for statim samples and at no cost for routine samples.
Finally, as expected, the fixed-schedule policy yields the best results in terms of the median patient TAT for statim and routine samples.

Compared with the threshold-based policy (which is often used in laboratories), the improvements for vitals are even more pronounced---9.7 minutes for vital samples at the 0.95 quantile level and 4.9 minutes at the median patient TAT.
Although the benefit of using stochastic decision-making by the optimization MIQP model in terms of absolute numbers may not seem very large over the reactive look-ahead policy, which does not use distributional information, it needs to be considered that the centrifuge timing is the only degree of freedom in the considered problem, as depicted in Figure~\ref{fig:process-timing}.
Thus, the results should be interpreted in terms of relative patient TAT reduction with respect to a cycle time of $\delta=15\cdot60$ seconds.
In this view, the stochastic MIQP policy virtually ``shortens'' the cycle time of the centrifuge by 65\% compared with the currently deployed threshold-based policy or 25\% compared with the look-ahead policy.

We also evaluated the achieved laboratory TAT to assess how it differs from the patient TAT.
When we compute the patient-laboratory TAT difference for vital samples for different algorithms at the 0.95 quantile level, it can be seen that this gap is not constant and is strictly increasing with decreasing amounts of information exposed to the algorithm.
Furthermore, when focusing on the median laboratory TAT for vital samples, the stochastic MIQP and look-ahead policies achieve the same laboratory TAT but different patient TAT values.
Therefore, these observations suggest that patient and laboratory TATs, although correlated, are generally different quantities, and the transport time of samples should not be overlooked.

Finally, we discuss the computational time consumed by the solution of the MIQP model.
We impose a time limit of 10 seconds for solving the model as a safeguard to ensure its real-time performance.
However, the Gurobi solver in version 11 presents substantial improvements in the solution of nonconvex MIQPs, and on the considered dataset, it typically achieves run times below 1~second.
Therefore, its efficiency for the number of samples considered in this case study is sufficient.
If the number of samples increases over time, the model's performance could be further improved by providing a warm-started solution, since finding a feasible solution to the problem is easy.
However, we have not explored this option further, as performance is not an issue.

\subsection{Different objective functions}\label{sec:exp_objectives}
In this experiment, we are interested in how different objective functions applied inside the stochastic MIQP model introduced in Section~\ref{sec:mip_alg} influence the patient TAT distribution.
Specifically, we investigate the use of the tardiness function and the related objectives, such as the total squared tardiness and the number of tardy jobs, and their impact on the patient TAT of the vital samples.

To introduce the concept of tardy samples, we need to define so-called virtual due dates.
The virtual due date for a sample $j\in\mathcal{S}^{\textsc{vital}}$ is the maximum patient TAT that we allow sample $j$ to reach without incurring any loss.
Thus, it is given as the time from its registration in the system $\hat{r}_j$ plus a given slack time of $\beta$ minutes that the user is willing to accept as a completion time.
Having virtual due dates defined, one can work with objective functions, such as the total tardiness $\sum_j T_j = \sum_j \max\{0, C_j -  \hat{r}_j-\beta\}$, squared tardiness $\sum_j T_j^2$, or the number of tardy jobs $\sum_j U_j = \sum_j \mathbb{1}(C_j -  \hat{r}_j-\beta)$.
\begin{figure}[ht]
    \centering
    \pgfplotsset{
        cycle list/PuBu-8,
        cycle multiindex* list={
            mark list*\nextlist
            PuBu-8\nextlist
        },
    }
   
    \begin{tikzpicture}
        \begin{axis}[
            height=0.3\textwidth, width=0.6\textwidth,
            legend pos={outer north east},
            legend cell align=left,
            xmajorgrids, ymajorgrids,
            xmax=65,
            xlabel={Patient TAT [min]}, ylabel={vital CDF [-]},legend columns=2, 
            label style={font=\scriptsize},
            tick label style={font=\scriptsize},
         colormap name=greenyellow,
            cycle multiindex* list={
                [samples of colormap=12]\nextlist
                mark list\nextlist
            },
            ]
            \addplot[thick,red,mark=none] table [x=probability, y=sota_vital_priority_with_timeout_algorithm, col sep=comma]{s53_nonhierarchical_ecdf.txt};
            \addlegendentryexpanded{look-ahead policy}
            \addlegendimage{empty legend}\addlegendentry{}
            \addlegendimage{empty legend}
            \addlegendentry{\hspace{-.7cm}{stochastic MIQP $\sum T_j$:}}
            \addlegendimage{empty legend}\addlegendentry{}
            \foreach \k in {40,...,50}{
                \edef\tempcol{miqp_cond_nonhiear_dd\k}
                \addplot+[thick,solid,mark=none] table [x=probability, y=\tempcol, col sep=comma]{s53_nonhierarchical_ecdf.txt};
                \addlegendentryexpanded{$\beta=\k$ mins}
            }
            \addplot+[thick,solid,mark=none,blue!80] table [x=probability, y=miqp_cond_hiear_dd0, col sep=comma]{s53_nonhierarchical_ecdf.txt};
            \addlegendentryexpanded{$\beta=0$ mins}
        \end{axis}
    \end{tikzpicture}
    \caption{The performance for vital samples of the model with the total tardiness objective.}
    \label{fig:exp_tardiness}
\end{figure}
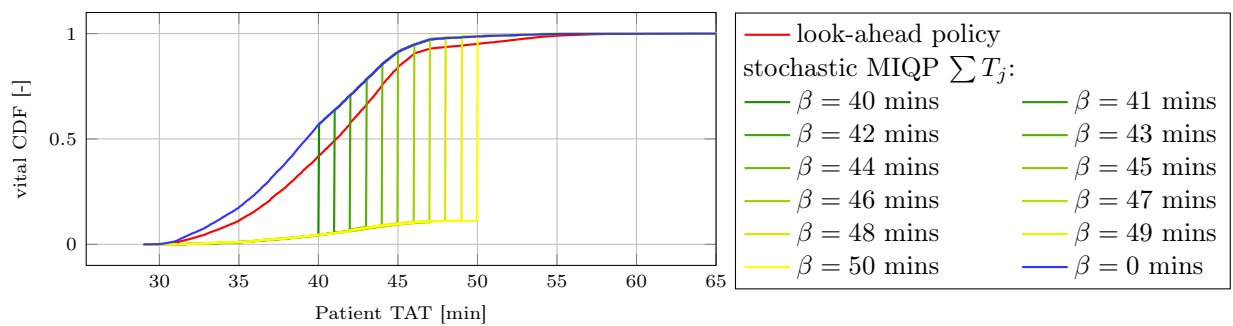

The motivation behind the total tardiness objective is to increase the decision space of the algorithm by allowing it to postpone the centrifuge run to anticipate the registration and transport of new samples.
The results are shown in Figure~\ref{fig:exp_tardiness} for values of virtual due date $\beta\in\{0, 40, 41, \ldots, 50\}$ minutes.
Furthermore, note that the total tardiness objective with a value of $\beta=0$ minutes resembles total flow time minimization since $\hat{r}_j \leq C_j$.
A closer inspection of the results suggests that the distribution of patient TAT becomes narrower since the model exploits the tardiness function and runs the centrifuge so that most samples are completed exactly $\beta$ minutes after their registration.
The narrower distribution of TATs is generally better, as laboratories and their clients (i.e., medical doctors) prefer a more predictable TAT.
Therefore, it would be meaningful to optimize, e.g., the difference between the 0.75 and 0.25 quantiles of the patient TAT distribution, to improve patient TAT consistency. 
However, we do not further explore this possibility in this paper.

Similar outcomes were observed for the sum of the squared tardiness function. 
%
The pronounced effect on the distribution of patient TATs seems to reflect the (expected) number of tardy jobs criterion. In this case, the solution is more sensitive to the actual value of the virtual due date $\beta$; however, it is still superior to the total flow time minimization.
From the observations above, it appears that the stochastic MIQP model achieves similar results as long as the objective function is strictly nondecreasing (i.e., $f(C_1, \ldots, C_n)$ is unbounded as $C_j \xrightarrow{} +\infty$) in the completion times $C_j$ of the samples.

\subsection{Perfect-knowledge offline setting}\label{sec:ex_pk_offline}
In this experiment, we compute some upper bounds on the performance that any online policy can achieve.
To assess this, we use the perfect-knowledge offline algorithm proposed in Section~\ref{sec:perfect_knowledge}.
This algorithm has access to all realizations of registration $\hat{r}_j$ and transportation times $\hat{\tau}_j$ in advance and thus makes globally optimal decisions concerning patient TAT minimization for vital samples and---with some assumptions discussed further---samples of other priorities.
To make the computations tractable, we give the perfect-knowledge algorithm information about the samples to be released over a one-day horizon.
Therefore, the achievable patient TATs are optimized for every day separately.
However, in the real data, we observed a gap of several hours between the last sample of a day and the first sample of the following day.
Considering the centrifuge cycle time of 15 minutes, this gap effectively voids the dependencies between the days and makes the laboratory process memoryless on the scale of the 24-hour horizon.
Thus, we conclude that solving each day separately is, in practice, essentially identical to the solution that considers the whole horizon worth months of data.

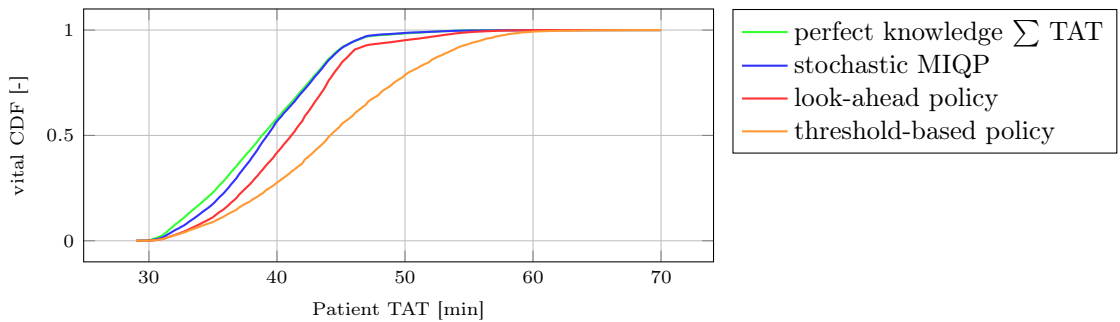
\begin{figure}[ht]
    \centering
    \pgfplotsset{
        cycle list/PuBu-8,
        cycle multiindex* list={
            mark list*\nextlist
            PuBu-8\nextlist
        },
    }
 
    \begin{tikzpicture}
        \begin{axis}[
            height=0.3\textwidth,
            width=0.6\textwidth,
            legend pos={outer north east},
            legend cell align=left,
            xmajorgrids, ymajorgrids,
            xlabel={Patient TAT [min]},
            ylabel={vital CDF [-]},
            label style={font=\scriptsize},
            tick label style={font=\scriptsize},
            table/x=probability,
            cycle list={green!80, blue!80, red!80, orange!80, cyan!80}
            ]
            \foreach \y in {sum_max_max,miqp_cond_hiear_dd0,sota_vital_priority_with_timeout_algorithm,sota_algorithm}{
                \addplot+[thick,mark=none] table [y=\y, col sep=comma]{s54_result.txt};
            }
            \addlegendentry{perfect knowledge $\sum$ TAT}\addlegendentry{stochastic MIQP}\addlegendentry{look-ahead policy}\addlegendentry{threshold-based policy}\addlegendentry{fixed-schedule policy}
        \end{axis}
    \end{tikzpicture}
    \caption{Patient TAT for vital samples obtained by the offline perfect-knowledge algorithm compared to online policies without perfect knowledge of all parameters.}\label{fig:offline_vs_online}
\end{figure}

We use 100 months' worth of data to estimate this upper bound on the performance of any online policy on the vital samples. 
The perfect-knowledge offline algorithm, which minimizes the daily sum of TATs, has access to all the realizations in advance and thus possesses complete information.
The fixed-schedule, threshold-based, look-ahead, and stochastic MIQP policies are executed in the online setting without access to all realizations in advance, which is the same as in the previous sections.
The distribution of patient TATs of vital samples resulting from this experiment is shown in Figure~\ref{fig:offline_vs_online}.
The most important outcome is that with perfect knowledge of all samples, slightly better patient TAT results are attainable for low quantiles. 
Nevertheless, for the high quantiles, the results of the offline algorithm with perfect knowledge almost match the performance of our online stochastic MIQP policy, which does not have access to future realizations of registration and transportation times.
This perhaps surprising result also suggests that, for example, even an optimal predictor of transportation times would not significantly improve the performance of the online stochastic MIQP policy.

Computing an upper bound on the performance of the algorithms for other sample priorities, i.e., \textsc{statim} and \textsc{routine}, is more challenging for two reasons: (i)~the problem resembles a hierarchical optimization problem; i.e., the solution adopted for \textsc{vital} samples affects the solution space of \textsc{statim} samples and so on for \textsc{routine}, and (ii)~a large number of \textsc{statim} and \textsc{routine} samples rules out the use of objective functions other than the minimization of the maximum TAT, as discussed in Section~\ref{sec:perfect_knowledge}.
As these aspects make the evaluation of the algorithms more complex, we provide more insights into the results and present them in two tables---Table~\ref{tab:ex_ub_statim} for statim and Table~\ref{tab:ex_ub_routine} for routine samples.
To increase the statistical significance of the results produced by the hierarchical optimization procedure, we used 100 months' worth of data in this experiment.

In Tables~\ref{tab:ex_ub_statim} and~\ref{tab:ex_ub_routine} below, each column in \textit{All offline} and \textit{All online} sections represents a sequence of algorithms that were used in the corresponding stages.
The algorithms we consider are offline with perfect knowledge (denoted as \textit{All offline}) or online policies (threshold-based or stochastic MIQP), denoted as \textit{All online}.
We compute statistics for the samples in a hierarchical manner.
That is, first, we compute statistics ($\max$, 0.95 quantile, and mean) over the samples in each day individually and then compute $\max$ and mean statistics over the aggregated days, as given in the \textit{Aggregation} column.
In this way, we report six numbers for each sequence of algorithms used in different stages (sample priorities) in the optimization.
The best results in each row are highlighted in bold.

The results in Table~\ref{tab:ex_ub_statim} show that the best achievable performance in terms of the maximum patient TAT for the statim samples is 244 minutes, which is only 7.2 minutes shorter than that of the stochastic MIQP, which has incomplete knowledge.
However, owing to the objective function focusing on the worst-case patient TAT, the performance of the perfect knowledge offline algorithms is less dominant in the 0.95 quantile measure. The outcomes confirm that the threshold-based policy is competitive for the statim samples, especially when it is optimized for the expected case.

With respect to the performance for routine samples in terms of the maximum TAT, the difference between the best offline perfect-knowledge algorithm and stochastic MIQP is only 4.4 minutes. 
Furthermore, the impact of using different online methods (threshold-based and stochastic MIQP) for routine samples alone is marginal.
This is because the starting times of the majority of batches are predetermined before routine samples because of the large number of statim samples optimized in the previous stage, which leaves little room for optimizing routine samples.

In any case, looking at the minimum along each row of Tables~\ref{tab:ex_ub_routine} and~\ref{tab:ex_ub_statim}, the numbers provide some guidance on what can be achieved with different sequences of full-knowledge offline algorithms and online policies combined with different aggregations of results.
From this, we estimate that our best online policy, the stochastic MIQP, delays routine samples by up to 7.2 minutes compared with the theoretical optimum. 
Interestingly, the mean TATs are similar for all the considered strategies and for both the statim and routine samples.
This suggests that optimization for different robust objectives, such as the daily maximum or 95th percentile, can be performed without a detrimental impact on the average case.
\begin{table}[ht]
    \centering
    \caption{Patient TATs in minutes of statim samples achieved by different sequences of algorithms.}
        \label{tab:ex_ub_statim}

     \resizebox{.85\textwidth}{!}{%
    \begin{tabularx}{\textwidth}{*{2}{l}*{2}{l}*{5}{L}}
        \toprule
        \multicolumn{2}{l}{Aggregation} & Stage & \multicolumn{2}{l}{All offline [min]} & \multicolumn{2}{l}{All online [min]} \\
        \cmidrule(lr){1-2} \cmidrule(lr){3-3} \cmidrule(lr){4-5} \cmidrule(lr){6-7} 
         Daily & Total & vital & $\max$~TAT & $\sum$ TAT & threshold-based & stochastic MIQP \\
         & & statim &$\max$~TAT & $\max$~TAT& threshold-based & stochastic MIQP \\
        \midrule
    max  & max &  & \textbf{244.0}     &   \textbf{244.0}     &      247.3      &          251.2 \\
         & mean & & \textbf{166.1}     &   \textbf{166.1}     &      171.9      &          172.2 \\[.25em]
    \ppp & max &  & \textbf{170.9 }    &   \textbf{170.9}     &      173.3      &          173.3 \\
         & mean & & \textbf{112.0}     &   112.1     &      112.2      &          112.4 \\[.25em]
    mean & max &  &  83.5     &    83.5     &       83.6      &           83.6 \\
         & mean & &  73.9     &    73.9     &       \textbf{73.5 }     &           73.7 \\
        \midrule
        \multicolumn{7}{l}{\scriptsize{Legend: $\sum$ TAT: Sec.~\ref{subsub:exactMILP}, $\max$ TAT: Sec.~\ref{subsub:maxTAT}, threshold-based: Sec.~\ref{sec:sota_alg}, stochastic MIQP:  Sec.~\ref{sec:mip_alg}.}}\\
       \bottomrule
    \end{tabularx}%
    }
\end{table}
\begin{table}[ht]
    \centering
    \caption{Patient TATs in minutes of routine samples achieved by different sequences of algorithms.}
        \label{tab:ex_ub_routine}
           \resizebox{.85\textwidth}{!}{%
    \begin{tabularx}{\textwidth}{*{3}{l}*{3}{L}}
        \toprule
        \multicolumn{2}{l}{Aggregation} & Stage & \multicolumn{1}{l}{All offline [min]} & \multicolumn{2}{l}{All online [min]} \\
        \cmidrule(lr){1-2} \cmidrule(lr){3-3} \cmidrule(lr){4-4} \cmidrule(lr){5-6} 
         Daily & Total & vital & $\max$~TAT & threshold-based & stochastic MIQP  \\
         & & statim & $\max$~TAT & threshold-based & stochastic MIQP  \\ 
         & & routine & $\max$~TAT  & threshold-based & stochastic MIQP \\
        \midrule
       max  & max  & & \textbf{259.4}        &     265.8  &     263.8    \\
            & mean & & \textbf{161.9}        &     167.1  &     167.3    \\[.25em]
       \ppp & max	& & \textbf{192.2}        &     196.2  &     196.2    \\
            & mean & & \textbf{123.7}        &     \textbf{123.7}  &     123.9    \\[.25em]
       mean & max	& & \textbf{133.8 }       &     144.6  &     144.6    \\
            & mean & &  \textbf{89.3}        &      89.4  &      89.7    \\
        \midrule
        \multicolumn{6}{l}{\scriptsize{Legend: $\max$~TAT: Sec.~\ref{subsub:maxTAT}, threshold-based: Sec.~\ref{sec:sota_alg}, stochastic MIQP: Sec.~\ref{sec:mip_alg}.}}\\
       \bottomrule
    \end{tabularx}%
    }
\end{table}
\subsection{Influence of system parameters}\label{sec:system_params}
The cycle time of the centrifuge is among the system parameters that significantly influence the patient TAT.
The cycle time of the centrifuge (including the loading and unloading of samples) depends on the manufacturer, its setup, and its capacity.
When comparing the quality of different batching policies, it should be noted that the perceived efficiency of a policy depends on the values of the system parameters, such as the cycle time, that are applied.
Intuitively, if the cycle time is very long, e.g., two hours, then the performance of any online policy would essentially converge to the fixed-schedule policy that runs the centrifuge every two hours, since there would be a high probability of having some sample that missed the centrifuge start.
Additionally, if the centrifuge is very fast, e.g., 1~minute, then the potential advantage of any nontrivial method is upper-bounded by 1~minute since the samples can enter the centrifuge any minute, assuming that it has sufficient capacity.
In contrast, the efficiency of a batching policy becomes apparent when the cycle time falls into a range that is not too short and not too long.
\begin{figure}[ht]
    \centering
    \pgfplotsset{
        cycle list/PuBu-8,
        cycle multiindex* list={
            mark list*\nextlist
            PuBu-8\nextlist
        },
    }
 
    \begin{tikzpicture}
        \begin{axis}[
            height=0.3\textwidth, width=0.6\textwidth,
            legend pos={outer north east}, legend cell align=left,
            xmajorgrids, ymajorgrids,
            xlabel={Patient TAT [min]}, ylabel={vital CDF [-]},
            label style={font=\scriptsize}, tick label style={font=\scriptsize},
            table/x=probability,
            ]
            \addlegendimage{empty legend}
            \addlegendentry{\hspace{-.7cm}{stochastic MIQP:}}
            \addplot+[thick,mark=none, blue!80] table [y=20, col sep=comma]{s55_result_miqp_cond_hiear_dd0.txt};
            \addplot+[thick,mark=none, blue!60] table [y=17, col sep=comma]{s55_result_miqp_cond_hiear_dd0.txt};
            \addplot+[thick,mark=none, blue!40] table [y=14, col sep=comma]{s55_result_miqp_cond_hiear_dd0.txt};
            \addplot+[thick,mark=none, blue!20] table [y=10, col sep=comma]{s55_result_miqp_cond_hiear_dd0.txt};
            \addlegendentry{$\delta=20$ mins}\addlegendentry{$\delta=17$ mins}\addlegendentry{$\delta=14$ mins}\addlegendentry{$\delta=10$ mins}
            \addlegendimage{empty legend}
            \addlegendentry{\hspace{-.7cm}{look-ahead:}}
            \addplot+[thick,mark=none, red!80] table [y=20, col sep=comma]{s55_result_timeout_algorithm.txt};
            \addplot+[thick,mark=none, red!60] table [y=17, col sep=comma]{s55_result_timeout_algorithm.txt};
            \addplot+[thick,mark=none, red!40] table [y=14, col sep=comma]{s55_result_timeout_algorithm.txt};
            \addplot+[thick,mark=none, red!20] table [y=10, col sep=comma]{s55_result_timeout_algorithm.txt};
            \addlegendentry{$\delta=20$ mins}\addlegendentry{$\delta=17$ mins}\addlegendentry{$\delta=14$ mins}\addlegendentry{$\delta=10$ mins}
        \end{axis}
    \end{tikzpicture}
    \caption{Patient TAT distribution of vital samples for different centrifuge cycle times~$\delta$.}
    \label{fig:ex_spinnig_time}
\end{figure}
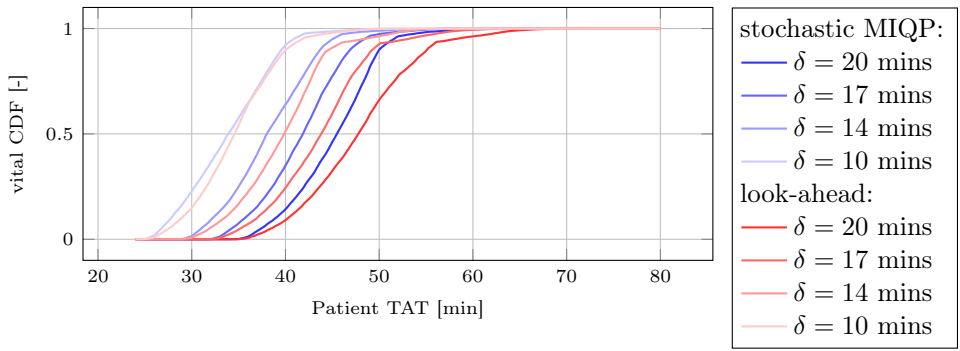

\begin{figure}[ht]
    \centering
        \begin{tikzpicture}
        \begin{axis}[
            height=0.28\textwidth,
            width=0.6\textwidth,
            xmajorgrids, ymajorgrids,
            xlabel={centrifuge cycle time $\delta$ [min]},
            ylabel={Patient TAT [min]},
            legend pos={outer north east},
            legend cell align=left,
            label style={font=\scriptsize},
            tick label style={font=\scriptsize},
            legend style={cells={align=left}}
            ]
            \addplot[smooth,red,thick] table [x=DELTA, y=QPatientTAT, col sep=comma] {s55_summary_timeout_algorithm.txt};
            \addlegendentry{look-ahead policy,\\ 0.95 quantile}

            \addplot[smooth,blue,thick] table [x=DELTA, y=QPatientTAT, col sep=comma] {s55_summary_result_miqp_cond_hiear_dd0.txt};
            \addlegendentry{stochastic MIQP,\\ 0.95 quantile}
            
            \addplot[smooth,red,dashed,thick] table [x=DELTA, y=MedianPatientTAT, col sep=comma] {s55_summary_timeout_algorithm.txt};
            \addlegendentry{look-ahead policy,\\ 0.5 quantile}

            \addplot[smooth,blue,dashed,thick] table [x=DELTA, y=MedianPatientTAT, col sep=comma] {s55_summary_result_miqp_cond_hiear_dd0.txt};
            \addlegendentry{stochastic MIQP,\\ 0.5 quantile}
            
        \end{axis}
    \end{tikzpicture}
    \caption{The 0.95 and 0.5 quantiles of patient TAT for different centrifuge cycle times.}
    \label{fig:ex_spinning_95q}
\end{figure}
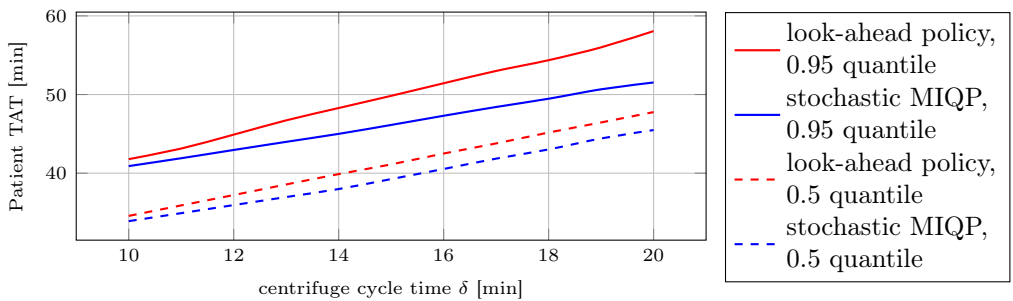

To evaluate the efficiency of different alternative system configurations, we created separate simulations with centrifuge cycle times of $\delta\in\{10,14,17,20\}$ minutes.
A comparison between the look-ahead policy and the stochastic MIQP model is shown in Figure~\ref{fig:ex_spinnig_time}.
The two methods perform similarly well in terms of the best cases and lower quantiles of the patient TAT.
However, with increasing cycle time, stochastic MIQP gains greater advantages in terms of higher quantiles as the precise timings of the centrifuge become more critical to avoid penalties for sample miss. 
This is apparent in the graph in Figure~\ref{fig:ex_spinning_95q}, which depicts the values of the 0.95 and 0.5 quantiles of the patient TAT depending on the value of $\delta\in\{10, 11, \ldots, 20\}$.
With respect to the 0.95 quantile, the results suggest that with increasing cycle time, stochastic MIQP gains more advantages, whereas for shorter cycle times below 10 minutes, the advantage of stochastic MIQP essentially diminishes.
The difference in the median performance is largely preserved across different cycle times.

\subsection{Summary}\label{sec:insights}
One of the main outcomes of this study is the assessment of how different information contexts can improve the patient TAT of vital samples.
In that regard, these results are best seen in Figure~\ref{fig:offline_vs_online} and in Table~\ref{tab:ex_numerical_all}.
The difference between the curves corresponding to the threshold-based and look-ahead policies quantifies the benefits of using the information regarding sample registration in the system.
It can be seen that implementing a system that considers such information improves the 0.95 quantile of patient TAT by as much as 6.7 minutes (see Table~\ref{tab:ex_numerical_all}).
An additional 3.7-minute improvement can be achieved by utilizing the distributional information of transportation times and optimization-based decision-making, which avoids various rule-based logic and timeout policies in favor of an optimization model that continually reevaluates its decision on the basis of newly revealed information to determine the optimal centrifuge start time. 
These improvements are obtained at virtually no cost in terms of degradation of the patient TAT of statim samples (0.2 minutes, both for the 0.95 and 0.5 quantiles).
The perfect-knowledge MIQP curve shows an upper bound on the performance of any online policy with respect to vital samples.
Knowing the exact transportation times and registration times in advance cannot be used to  eliminate  the upper tail of patient TAT cases, but such information might be used to slightly improve patient TAT in the bottom 50\% of vital samples.
\section{Managerial insights}
Our research contributes to the management of laboratories and hospitals aiming at total laboratory automation (TLA), where the ultimate goal is to process samples without any assistance from a human operator except for clinical validation of the results.
Current automation systems typically distinguish between two sample priorities, i.e., statim and routine. However, if the hospital defines vital priorities, such samples are processed manually. In practice, this means that a laboratory employee takes the sample and manually performs all the operations with the sample (e.g., using manual input for individual analyzers). If multiple vital samples arrive at the laboratory simultaneously, the laboratory should have sufficient staff to process them. Therefore, the manual processing of vital samples increases the demands on human resources.
On the other hand, the manual processing of vital samples contradicts the TLA. Therefore, our research addresses the processing of vital samples and examines how they should be handled by laboratory automation--specifically, by a centrifuge. From a laboratory management point of view, this is an important step in eliminating manual and routine work (which is often a source of mistakes, such as incorrect blood in tube errors \cite{Rosenbaum2018}). With TLA, which also covers vital samples, highly skilled laboratory employees can be reassigned to more qualified positions, supporting the more efficient use of human resources.

This paper also concerns sample flow quality control. Hospitals often see laboratory TAT as the main key performance indicator for laboratory samples. 
The reason for this is that responsibility for the flow of samples--from the moment the sample is drawn until the result is reported--is distributed among multiple stakeholders. No single party has overall responsibility for or control over the entire process.
Individual hospital wards are responsible for collecting the samples and passing them to the transport. A technical department is responsible for transporting the samples, often via tube mail. 
Finally, the laboratory receives the samples, analyzes them, and reports the results back to the doctor. However, hospitals typically do not have a body that monitors and manages the entire process.
This aspect may result in a significant variation in the time it takes to deliver laboratory results, which may be perceived negatively by doctors. 
This opinion is supported by data from University Hospital Královské Vinohrady, presented in Figure~\ref{fig:transport-distribution}, which shows how transport times from different wards can be diverse; thus, one department may consider a certain laboratory TAT to be sufficient, while another may not.
Given that laboratory TATs may not provide all relevant information about the quality of care, we want to draw the attention of hospital authorities to patient TATs.
Hospitals should periodically monitor patient TAT and analyze how individual departments contribute to this time (see the timing diagram in Figure~\ref{fig:process-timing}). 
Moreover, hospitals should have a single entity responsible for managing the quality of the process and ensuring that the entire workflow is centrally monitored. 
Managers should also look for ways to automate the process and eliminate the need for human interaction with samples, which often causes delays and introduces unpredictability.
Of course, transport time may vary across different departments due to their distance from the laboratory. 
Therefore, our paper shows how information related to the transport of samples can be incorporated into laboratory automation systems to improve patient TAT, and it can serve as further motivation for automated sample processing. 
This should be in the economic interest of hospitals, as reducing patient TAT at high quantile levels contributes to shorter patient length of stay~\cite{holland2005reducing} and thus has an eminent impact on the quality of care provided by the hospital and healthcare costs.

\section{Conclusion and future research}
This paper addresses an important topic concerning the speed of delivery of laboratory test results. We highlight the current limitations of laboratories that evaluate laboratory TAT, not patient TAT. Therefore, we propose a new batching strategy for centrifuging samples that incorporates knowledge of transport times from individual wards. Experimental results on real data show that the proper batching of samples obtained by our strategy can significantly reduce the patient TAT of vital samples, i.e., samples of patients in immediate danger of death, while not degrading the patient TAT of less urgent samples.

Regarding the method proposed in this paper, several natural extensions of our work are possible.
First, we do not address the analytic part of the system.
If the load of analyzers is high, then a sample may need to wait before it is processed, and this waiting time depends on the composition of the sample methods to be evaluated within a certain short-term time window.
Therefore, it would be interesting to model the analytic part of the system, consider it in the patient TAT calculation, and utilize it by the batching strategy in a way that leads to balanced use of the analyzers.
Moreover, since handling vital samples is time-critical and the historical data may be insufficient to construct a complex model of the transport time distribution, a distributionally robust approach~\cite{NOVAK2022438} could be employed. Such methods are already applied in the medical field, for example, for surgical block allocation~\cite{wang2019distributionally}.
Next, the stochastic MIQP model could be extended to handle samples with different priorities, leading to optimized batching of non-vital samples as well. 
Finally, it would be interesting to extend the method to consider more than a single centrifuge.

Further research should also address the use of priorities in total laboratory automation. As described in the previous section, vital samples are often handled manually, contrary to the concept of total laboratory automation. Nevertheless, this is not the only issue that current laboratories face regarding sample priorities. Statim priorities, which should be used for patients with slightly less urgent needs than those requiring a vital indication, are often overused for samples that need to be processed, e.g., before a ward round. Then, the laboratory may first process samples from inpatients who are waiting for the ward round, and the sample of a patient taken to the hospital by ambulance and suffering from acute myocardial infarction (whose sample was also labeled statim) may be processed subsequently. In such a situation, care for the patient with acute myocardial infarction is delayed since the laboratory is not able to distinguish which statim sample is most urgent.
To the best of our knowledge, no existing work addresses total laboratory automation and the design of sample priorities.

\section*{Acknowledgments}
This work was supported by the Grant Agency of the Czech Republic under the Project \mbox{GACR~22-31670S} and by the European Union under the project ROBOPROX -- Robotics and Advanced Industrial Production (reg. no. CZ.02.01.01/00/22\_008/0004590).

\bibliography{sn-bibliography}
\begin{appendices}
\section{Perfect-knowledge optimization subroutines}
\label{app:offline}
The appendix provides a detailed description of the perfect-knowledge algorithm subroutines and an illustrative numerical example of the algorithm for solving a small instance. 
\subsection{Total TAT minimization subroutine}
\label{subsub:exactMILP}
We begin by describing a subroutine, which is the exact MIP formulation to be solved by general-purpose solver software.
The formulations can be used to optimize different objective functions if the number of samples is small enough. This section specifically discusses $\sum \textrm{TAT}_j$ minimization, i.e.,\ the total TAT, but can be easily adjusted for related objective functions.

To minimize the considered $\sum \textrm{TAT}_j$ objective, we describe the batching model by equations~\eqref{eq:milp:1}--\eqref{eq:milp:n}:
{
\medmuskip=0mu
\thinmuskip=-1mu
\thickmuskip=0.5mu
\nulldelimiterspace=0pt
\scriptspace=.5pt
\begin{alignat}{4}
    \min  \textstyle\sum_{j\in  \mathcal{J}} C_j&               &\ & \label{eq:milp:1}\\
    \shortintertext{subject to}
    \hat{r}_j + \hat{\tau}_j + \delta + p_j     &\leq C_j       & & \forall j \in  \mathcal{J} \\
    S^B_{b-1} + \delta                          &\leq S^B_{b}   & & \forall b \in \mathcal{B} \setminus \{1\}\\
    \textstyle\sum_{b\in \mathcal{B}} y_{j, b}  &= 1            & & \forall j \in  \mathcal{J}\\
    \textstyle\sum_{j\in  \mathcal{J}} y_{j, b} &\leq N         & & \forall b \in \mathcal{B}\\
    S^B_b  + \delta + p_j + M(1-y_{j,b})        &\geq C_j       & & \forall b\in \mathcal{B}, \ \forall j \in  \mathcal{J}\label{eq:offlineMilp:c1}\\
    S^B_b  + \delta + p_j - M(1-y_{j,b})        &\leq C_j       & & \forall b\in \mathcal{B}, \ \forall j \in  \mathcal{J}\label{eq:offlineMilp:c2}\\
    \shortintertext{where}
    &&\mathclap{y \in \{0, 1\}^{|\mathcal{J}| \times |\mathcal{B}|}, \ C \in \mathbb{R}^{| \mathcal{J}|},\ S^B \in \mathbb{R}^{|\mathcal{B}|},\qquad\qquad\qquad\qquad\qquad } \label{eq:milp:n}
\end{alignat}}%
where $M$ is a big-M constant and $\mathcal{J}$ is a set of samples to be scheduled by the subroutine. We denote by $\mathcal{B}$ the set of sample batches (i.e., centrifuge runs) considered. We upper bound the number of batches as the maximum number of batches in a day (based on the batch time) or by the number of samples; hence, the set is given by 
\begin{equation*}
    \mathcal{B} = \big\{1, 2, \dotsc, \min\{|\mathcal{J}|,\ \lfloor 24 \cdot 60 \cdot 60 / \delta \rfloor\}\big\}.
\end{equation*}
Binary variables $y_{j, b}$ determine whether sample $j\in  \mathcal{J}$ is included in batch $b\in \mathcal{B}$.
The variables $S^B_b$ denote the start time of batch $b$.
The objective function is to minimize the sum of completion times, which is equivalent to the minimization of the sum of patient TATs.
After solving the model, we can derive batch start times and sample-batch assignments from the variables $y_{j, b}$ and $S^B_b$ for $j\in\mathcal{J},\ b\in\mathcal{B}$.

The practical tractability of the model depends on the number of samples in $\mathcal{J}$. 
In our case study, the approach could be used only in the first stage since the daily number of statim samples is almost 200, but the number of vital samples in our dataset is within the reach of the current MIP solvers. 
Furthermore, some modern MIP solvers can return all possible optimal solutions if their number is manageable. 
In our dataset, the number of unique optimal solutions for scheduling vital samples never exceeded 2, with over 99\% of days admitting only a single unique optimal schedule. 
Therefore, all of the solutions may be passed to later stages without affecting tractability, which provides a greater diversity of partial schedules to consider and may lead to improved overall results.

If the number of samples of other priorities were not large, the subroutine could also be used at different stages. 
To guarantee that the starting times of the samples from the previous stages are not changed in the new solution, one can simply freeze their completion times in equations~\eqref{eq:offlineMilp:c1} and~\eqref{eq:offlineMilp:c2}. Alternatively, one can introduce only constraints on the maximum start times, allowing the samples to be scheduled earlier. 
\subsection{Maximum TAT minimization subroutine}
\label{subsub:maxTAT}
In this section, we describe a polynomial-time algorithm for minimizing the maximal patient TAT. 
Essentially, it resembles a binary search on the maximal TAT value combined with a fast algorithm of \cite{Condotta2010} for checking feasibility. 
Owing to its favorable time complexity, it can be used to construct a subroutine for any stage of optimization.
The pseudocode is shown in Algorithm~\ref{alg:revised_tat_optimize}.
\begin{algorithm}
    \DontPrintSemicolon
    \caption{Subroutine for minimizing maximal TAT.}%
    \label{alg:revised_tat_optimize}
    
    \SetKwFunction{FExact}{Exact\_MILP}
    \SetKwFunction{FOptimize}{Optimize}
    \SetKwFunction{FFindFeasible}{Find\_Feasible\_Solution}
    \SetKwProg{Fn}{Function}{:}{\KwRet}
    \KwData{Set of considered samples $\mathcal{J}$, the set of scheduled samples $\mathcal{P}$ from previous stage, completion times from previous stages $C^{P}_i,\ \forall i\in\mathcal{P} \subseteq \mathcal{J}$, the vector of available times $\boldsymbol{\hat{a}}= (\hat{a}_1,\ldots, \hat{a}_n)$.}
    \KwResult{Start times of samples in $\mathcal{J}$.}
    
    \Fn{\FOptimize{$C^P$, $\mathcal{P}$, $\mathcal{J}$}}{
        \ForEach{$j \in \mathcal{J}$}{
            \uIf(\tcp*[h]{Scheduled in the previous stages}){$j \in \mathcal{P}$}{
                $d^\prime_j \gets C^P_j - p_j$\;
            }
            \Else(\tcp*[h]{Unscheduled samples are processed within 24h}){
                $d^\prime_j \gets \hat{a}_j + 24 \cdot 60 \cdot 60$\;
            }
        }
        $\boldsymbol{d}^{\text{copy}} \gets \boldsymbol{d^\prime}= (d^\prime_1,\ldots, d^\prime_{n})$\;
        $LB \gets 0,\quad UB \gets \max_{j \in \mathcal{J}}\{d^\prime_j - \hat{a}_j\}$\;
        \While{$LB \neq UB$}{
            $c \gets \lfloor(LB + UB) / 2\rfloor$\;
            \ForEach{$j \in \mathcal{J}$}{
                $d^\prime_j \gets \min\big\{ d^{\text{copy}}_j,\ \hat{a}_j + \delta + \max\{0,\ c - (p_j + \hat{\tau}_j)\}\big\}$\;
            }
            \eIf{\FFindFeasible{$\boldsymbol{\hat{a}}, \boldsymbol{d^\prime}$} exists}{
                $LB \gets c + 1$\;
            }{
                $UB \gets c$\;
            }
        }
        $c \gets LB$\;
        \ForEach{$j \in \mathcal{J}$}{
            $d^\prime_j \gets \min\big\{ d^{\text{copy}}_j,\ \hat{a}_j + \delta + \max\{0,\ c - (p_j + \hat{\tau}_j)\}\big\}$\;
        }
        \Return \FFindFeasible{$\boldsymbol{\hat{a}}, \boldsymbol{d^\prime}$}\;
    }
\end{algorithm}

First, deadlines are constructed for all the samples in $\mathcal{J}$. 
A deadline $d^\prime_j$ for $j \in \mathcal{J}$ indicates that sample $j$ cannot leave the centrifuge later than time $d^\prime_j$.
For the already scheduled samples $\mathcal{P}$ from the previous stage, the deadlines are constructed using the solution of the previous stages $C^P$. In this way, we guarantee that the schedule constructed by the current stage for all samples $\mathcal{J}$ must not postpone any sample from $\mathcal{P}$. 
For the remaining samples $j\in \mathcal{J}\setminus \mathcal{P}$, the initial deadlines are set to $d^\prime_j=\hat{a}_j + 24\cdot60\cdot60$, i.e., 24 hours after a sample's arrival at the laboratory. 
We assume that a feasible solution exists for such a set of deadlines since we consider instances where the laboratory can process all the daily requests (constituting an upper bound). 
To optimize the patient TAT for these samples, we also consider a second set of deadlines that are tighter but likely infeasible, e.g., $\hat{a}_j + \delta$ (see lines 12 and 19 of Algorithm~\ref{alg:revised_tat_optimize}).
These values represent lower bounds on the feasible deadlines, iteratively enlarged by the nonnegative term 
\begin{equation}\label{eq:cterm}
    \max\{0, c - (p_j + \hat{\tau}_j)\},
\end{equation}
whose value is selected using binary search to find a minimum~$c$ such that the problem remains feasible.

In the algorithm, $LB$ and $UB$ define bounds on $c$, tightened by the binary search until $LB=UB$, with the final value of $c$ derived.
To improve the quality of a solution, if sample $j$ reaches its minimal feasible deadline value, it is not further affected by~$c$, as given by the max term of Equation~\eqref{eq:cterm}. 
To account for postcentrifuging processing time and transport time, the terms $p_j$ and $\hat{\tau}_j$ are subtracted from $c$ in the same equation. 
In this way, samples with longer processing and transport times must have shorter deadlines.
Checking whether a solution for a given set of deadlines $\boldsymbol{d^\prime}=(d^\prime_1, \ldots, d^\prime_{n})$ can be constructed is done using \texttt{Find\_Feasible\_Solution}, an algorithm proposed in~\cite{Condotta2010} with $O(n^2)$ time complexity. 
The algorithm \texttt{Find\_Feasible\_Solution} returns a feasible solution (given by completion times of all samples) or reports that such a solution does not exist.   
The deadlines for the scheduled samples are minimized between lines 9--16 of Algorithm~\ref{alg:revised_tat_optimize}.

Below, we show why the algorithm optimizes the maximum patient TAT.
Consider completion times $C_j$ of samples $j\in \mathcal{J} \setminus \mathcal{P}$. 
If the solution exists, the deadlines $d^\prime_j$ guarantee that $C_j \leq d^\prime_j$. 
Assuming that $\omega_j$ is the time sample $j$ must wait before being centrifuged, we have $C_j = \hat{r}_j + \hat{\tau}_j + \omega_j + \delta$. With the knowledge that $\hat{a}_j = \hat{r}_j + \hat{\tau}_j$, the deadlines are computed in the algorithm by the formula
\begin{equation*}
    d^\prime_j =\min\big\{ d^{\text{copy}}_j,\ \hat{r}_j + \hat{\tau}_j + \delta + \max\{0,\ c - (p_j + \hat{\tau}_j)\}\big\}.
\end{equation*}
The outer minimum guarantees that the initial deadlines $d^{\text{copy}}_j$ cannot be extended, which preserves the solution from the previous stages. Therefore, for the scheduled samples, we can simply write
\begin{equation*}
    d^\prime_j = \hat{r}_j + \hat{\tau}_j + \delta + \max\{0,\ c - (p_j + \hat{\tau}_j)\}.
\end{equation*}
If $0 > c - (p_j + \hat{\tau}_j)$, the lower bound of the maximal patient TAT, $\max_j (\hat{\tau}_j + \delta + p_j)$, is attained. 
Otherwise, we have
\begin{equation*}
    d^\prime_j = \hat{r}_j + \hat{\tau}_j + \delta + (c - (p_j + \hat{\tau}_j)).
\end{equation*}
Then, for all $j$:
\begin{align*}
    C_j &\leq d^\prime_j \\
    \hat{r}_j + \hat{\tau}_j + \omega_j + \delta &\leq\hat{r}_j + \hat{\tau}_j + \delta + c - (p_j + \hat{\tau}_j) \\
    x_j + p_j + \hat{\tau}_j &\leq c \\
    \text{TAT}_j &\leq c.
\end{align*}
Since the algorithm terminates with the lowest value of $c$ for which a feasible solution exists, the lowest value of the maximal patient TAT is achieved. 
Finally, the solution with the minimized deadlines is constructed by \texttt{Find\_Feasible\_Solution}, and the start times of the batches are derived. The number of executions of \texttt{Find\_Feasible\_Solution} is pseudopolynomial and depends on the length of the scheduling horizon $UB$, e.g., in our case, $\log_2 UB \leq \log_2 (24 \cdot 60 \cdot 60 + \delta)$. 
Therefore, the overall time complexity of the algorithm is $O(n^2\cdot\log_2 UB)$.

The subroutine described in this section computes only a single optimal solution. 
However, the considered instances of the problem with maximum daily TAT minimization usually have a large number of optimal solutions, except for the vital stage.
Thus, lower priority stages of the offline perfect-knowledge algorithm (i.e., Algorithm~\ref{alg:revised_full_algorithm}) might be solved suboptimally.
The impact of multiple optimal solutions on the following stage is illustrated in the next section using an example. 
\subsection{Numerical example}
To demonstrate the functionality of the proposed algorithmic schema and its subroutines, we present a numerical example.
Consider an instance with 6 samples, whose parameters are summarized in Table~\ref{tab:fullinfo_example}. The instance is to be solved using the proposed algorithm schema, with $\sum \mathrm{TAT}_j$ minimization for the vital samples (using the subroutine from Sec.~\ref{subsub:exactMILP}) and maximum TAT minimization for the statim and routine samples (using the subroutine from Sec.~\ref{subsub:maxTAT}). 
First, the model \eqref{eq:milp:1}--\eqref{eq:milp:n} is used to schedule vital samples $\{1,2\}$. 
The results with two possible optimal solutions are displayed in Figure~\ref{fig:fullinfo_example}.
Solution 1 archives the sum of TATs equal to 60 minutes for the vital samples and schedules them in two separate batches.
As a result, $d^\prime_1 = 17.5$ and $d^\prime_2 = 32.5$. Then, the schedule of statim samples is constructed while minimizing the maximum patient TAT. 
The initial $d^\prime_j$ is set to $\hat{a}_j + UB$ (a large number) for $j\in\{3, 4, 5\}$. 
Subsequently, $c$~is minimized such that $c=2400$ (i.e., 40 minutes). The final values of $d^\prime_j$ (as computed at line 20 of Algorithm~\ref{alg:revised_tat_optimize}) can be found in Table~\ref{tab:fullinfo_example}. Note that sample $4$ is not centrifuged immediately since the centrifuge waits for sample $5$, resulting in a gap between the 2nd and 3rd batches, as shown in Figure~\ref{fig:fullinfo_example}. 
Finally, the routine sample is scheduled for the 4th centrifuge run with $c=4200$ (i.e., 70 minutes).
\begin{table}[h]
    \caption{An example instance and optimal solutions computed by Algorithm~\ref{alg:revised_full_algorithm}.}
        \label{tab:fullinfo_example}
    \centering
    \resizebox{.65\textwidth}{!}{
    \begin{tabular}{ccccccccccc}
        \toprule
        \multicolumn{6}{c}{Instance} & \multicolumn{2}{c}{Solution 1} & \multicolumn{2}{c}{Solution 2} \\
        \cmidrule(lr){1-6} \cmidrule(lr){7-8} \cmidrule(lr){9-10}
        $j$ & $\hat{r}_j$ & $\hat{\tau}_j$ & $\hat{a}_j$ & $p_j$ & $\lambda_j$ & $d^\prime_j$ & TAT & $d^\prime_j$ & TAT \\
            & [min]       & [min]          & [min] &  [min] & & [min] & [min] & [min] & [min] \\
        \midrule
        1 & 0.0  & 2.5  & 2.5  & 5 & \textsc{vital}    & 17.5 & 22.5 & 25.0 & 30.0 \\
        2 & 0.0  & 10.0 & 10.0 & 5 & \textsc{vital}    & 32.5 & 37.5 & 25.0 & 30.0 \\
        3 & 0.0  & 17.0 & 17.0 & 5 & \textsc{statim}   & 35.0 & 37.5 & 40.0 & 45.0 \\
        4 & 15.0 & 17.0 & 32.0 & 5 & \textsc{statim}   & 50.0 & 40.0 & 55.0 & 45.0 \\
        5 & 17.0 & 18.0 & 35.0 & 5 & \textsc{statim}   & 52.0 & 38.0 & 57.0 & 43.0 \\
        6 & 0.0  & 40.0 & 40.0 & 5 & \textsc{routine}  & 65.0 & 70.0 & 55.0 & 60.0 \\
        \bottomrule
    \end{tabular}%
    }
\end{table}
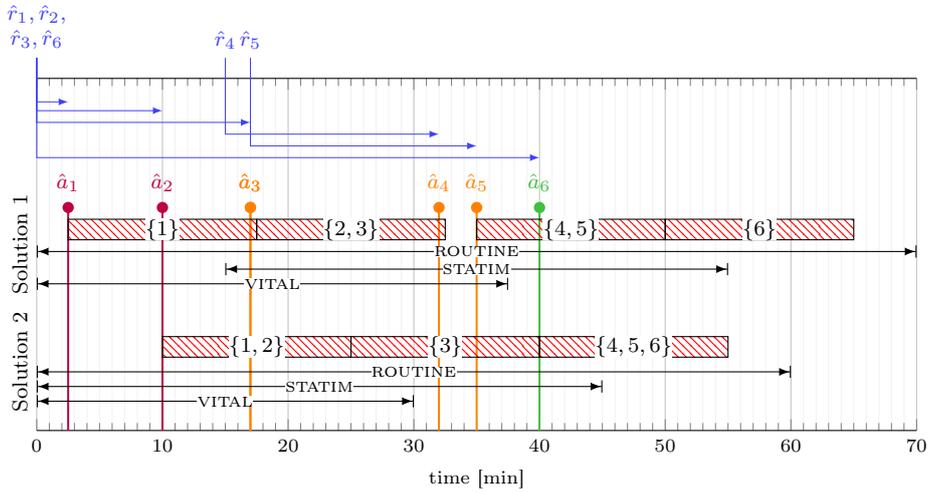
\begin{figure}[h]
    \centering
    \begin{tikzpicture}[
        every node/.style={font=\footnotesize},
        s1bar/.style={pattern=north west lines, pattern color=red},
        s2bar/.style={pattern=north west lines, pattern color=red},
        bartext/.style={pos=.5, font=\footnotesize, fill=white, inner sep=0},
        anode/.style={circle, fill, font=\footnotesize, inner sep=1.5pt},
        aline/.style={line width=.75pt},
        tatlabel/.style={midway, font=\footnotesize, fill=white, inner sep=0},
        tatline/.style={{| latex-latex[] |}},
        tauline/.style={line width=.25pt,-latex, blue!75},
        taulabel/.style={font=\footnotesize, text width=1cm, anchor=south, align=center},
    ]
        \begin{axis}[height=0.38\textwidth, width=0.8\textwidth,
            ymin=0, ymax=1.2, xmin=0, xmax=70,
            xlabel={time [min]},
            label style={font=\scriptsize},
            tick label style={font=\scriptsize},
            minor x tick num=9, hide y axis,
            grid=both,
            grid style={line width=.1pt, draw=gray!10},
            major grid style={line width=.2pt,draw=gray!50}, clip=false]

            \coordinate (alevel) at (axis cs: 0, 0.76);
            \draw[aline, purple] (axis cs: 2.5, 0) -- +(alevel) node[anode, label=$\hat{a}_1$] (rt1) {};
            \draw[aline, purple] (axis cs: 10,  0) -- +(alevel) node[anode, label=$\hat{a}_2$] (rt2) {};
            \draw[aline, orange] (axis cs: 17,  0) -- +(alevel) node[anode, label=$\hat{a}_3$] (rt3) {};
            \draw[aline, orange] (axis cs: 32,  0) -- +(alevel) node[anode, label=$\hat{a}_4$] (rt4) {};
            \draw[aline, orange] (axis cs: 35,  0) -- +(alevel) node[anode, label=$\hat{a}_5$] (rt5) {};
            \draw[aline, green!50!gray]   (axis cs: 40,  0) -- +(alevel) node[anode, label=$\hat{a}_6$] (rt6) {};

            \coordinate (taulevel) at (axis cs: 0, 1.27);
            \draw[tauline] (axis cs: 0, 0)  ++(taulevel) 
                node[taulabel] (r1) {$\hat{r}_1, \hat{r}_2,$\\ $\hat{r}_3, \hat{r}_6$} -- ++(axis cs: 0, -1*0.05-0.1) node (c) {} --  (c -| rt1);
            \draw[tauline] (axis cs: 0, 0)  ++(taulevel) node[taulabel] (r2) {} -- ++(axis cs: 0, -2*0.04-0.1) node (c) {} --  (c -| rt2);
            \draw[tauline] (axis cs: 0, 0)  ++(taulevel) node[taulabel] (r3) {} -- ++(axis cs: 0, -3*0.04-0.1) node (c) {} --  (c -| rt3);
            \draw[tauline] (axis cs: 15, 0) ++(taulevel) node[taulabel] (r4) {$\hat{r}_4$} -- ++(axis cs: 0, -4*0.04-0.1) node (c) {} --  (c -| rt4);
            \draw[tauline] (axis cs: 17, 0) ++(taulevel) node[taulabel] (r5) {$\hat{r}_5$} -- ++(axis cs: 0, -5*0.04-0.1) node (c) {} --  (c -| rt5);
            \draw[tauline] (axis cs: 0, 0)  ++(taulevel) node[taulabel] (r6) {} -- ++(axis cs: 0, -6*0.04-0.1) node (c) {} --  (c -| rt6);

            \draw[aline, orange]  (axis cs: 17,  0) -- +(alevel) node[anode, label={$\hat{a}_3$}] (rt3) {};


            \coordinate (bsize) at (axis cs:15,0.07);
            
            \coordinate (s1s) at (axis cs:0,.65);
            \draw[s1bar] (axis cs:2.5,0)  ++(s1s) node[](s1t){} rectangle ++(bsize) node[bartext] {$\{1\}$};
            \draw[s1bar] (axis cs:17.5,0) ++(s1s) rectangle ++(bsize) node[bartext] {$\{2, 3\}$};
            \draw[s1bar] (axis cs:35,0)   ++(s1s) rectangle ++(bsize) node[bartext] {$\{4, 5\}$};
            \draw[s1bar] (axis cs:50,0)   ++(s1s) rectangle ++(bsize) node[bartext] {$\{6\}$};

            \draw[tatline] (axis cs: 0,  -0.15)  ++(s1s) node[](s1b){} -- ++(axis cs: 37.5, 0) node[tatlabel]{$ \textsc{vital}$};
            \draw[tatline] (axis cs: 15, -0.10)  ++(s1s) -- ++(axis cs: 40,   0) node[tatlabel]{$ \textsc{statim}$};
            \draw[tatline] (axis cs: 0,  -0.04)  ++(s1s) -- ++(axis cs: 70,   0) node[tatlabel]{$\textsc{routine}$};

            \coordinate (s2s) at (axis cs:0,0.25);
            \draw[s2bar] (axis cs:10,0) ++(s2s)node[](s2t){} rectangle ++(bsize) node[bartext] {$\{1, 2\}$};
            \draw[s2bar] (axis cs:25,0) ++(s2s) rectangle ++(bsize) node[bartext] {$\{3\}$};
            \draw[s2bar] (axis cs:40,0) ++(s2s) rectangle ++(bsize) node[bartext] {$\{4, 5, 6\}$};

            \draw[tatline] (axis cs: 0, -0.15)  ++(s2s) node[](s2b){} -- ++(axis cs: 30, 0) node[tatlabel]{$\textsc{vital}$};
            \draw[tatline] (axis cs: 0, -0.1)   ++(s2s) -- ++(axis cs: 45, 0) node[tatlabel]{$\textsc{statim}$};
            \draw[tatline] (axis cs: 0, -0.05)  ++(s2s) -- ++(axis cs: 60, 0) node[tatlabel]{$\textsc{routine}$};

            \path (s2b) -- (s2t -| s2b) node[midway, rotate=90, anchor=south, align=center, xshift=5pt, font=\footnotesize]{Solution 2}; 
            \path (s1b) -- (s1t -| s1b) node[midway, rotate=90, anchor=south, align=center, xshift=5pt, font=\footnotesize]{Solution 1}; 
        \end{axis}
    \end{tikzpicture}
    \caption{Different solutions and their maximum TATs, for instance, in Table~\ref{tab:fullinfo_example}. Arrival times of vital samples are marked in purple, statim in orange, and routine in green. Release times are depicted in blue.}
    \label{fig:fullinfo_example}
\end{figure}

However, another optimal solution is possible for the vital samples with respect to the total patient TAT time criterion.
In Solution~2, the first centrifuge run is postponed to also contain sample $2$. Such a situation arises, for example, whenever two vital samples arrive exactly 450 seconds apart. 
This difference results in lowering the maximum patient TAT for routine samples from 70 to 60 minutes compared to Solution~1 but increases the maximum patient TAT for statim samples from 40 to 45 minutes. Furthermore, the difference in the maximum patient TAT of vital samples can be decreased from 37.5 to 30 minutes, while their sum remains 60 minutes.
\end{appendices}
\end{document}